\documentclass[10pt,a4paper,twoside]{amsart}
\usepackage{amsmath}
\usepackage{amsfonts}
\usepackage{amssymb}
\usepackage[applemac]{inputenc}
\usepackage[english]{babel}
\usepackage{verbatim}
\usepackage{graphicx}
\usepackage[all]{xy}
\xyoption{poly}
\usepackage{enumitem}
\usepackage{multirow}
\usepackage{longtable}
\usepackage{booktabs}

\usepackage{stmaryrd}

\usepackage[usenames,dvipsnames]{color}

\usepackage[obeyspaces]{url}

\usepackage[colorlinks = true, allcolors = cyan]{hyperref}

	\usepackage{amsthm}
	\usepackage{aliascnt}

	\theoremstyle{plain}

	\newtheorem{main}{Theorem}
	
	\newaliascnt{main_cor}{main}
	
	\aliascntresetthe{main_cor}

	\newaliascnt{prop}{thm}
	\newtheorem{prop}[prop]{Proposition}
	\aliascntresetthe{prop}
	\newaliascnt{lem}{thm}
	\newtheorem{lem}[lem]{Lemma}
	\aliascntresetthe{lem}

	\newaliascnt{cor}{thm}
	\newtheorem{cor}[cor]{Corollary}
	\aliascntresetthe{cor}
	
	\newaliascnt{conj}{thm}
	
	\aliascntresetthe{conj}

	\theoremstyle{definition}
	\newaliascnt{defn}{thm}
	\newtheorem{defn}[defn]{Definition}
	\aliascntresetthe{defn}

	\theoremstyle{remark}
	\newaliascnt{notation}{thm}
	\newtheorem{notation}[notation]{Notation}
	\aliascntresetthe{notation}
	\newaliascnt{rem}{thm}
	\newtheorem{rem}[rem]{Remark}
	\aliascntresetthe{rem}
	\newaliascnt{ass}{thm}
	\newtheorem{ass}[ass]{Assumption}
	\aliascntresetthe{ass}

	\newaliascnt{exmp}{thm}
	\newtheorem{exmp}[exmp]{Example}
	\aliascntresetthe{exmp}

\usepackage[capitalise]{cleveref}

	\Crefname{equation}{}{}
	\crefname{equation}{}{}
	\Crefname{thm}{Theorem}{Theorems}
	\Crefname{main}{Theorem}{Theorems}
	\Crefname{main_cor}{Corollary}{Corollaries}
	\Crefname{prop}{Proposition}{Propositions}
	\Crefname{lem}{Lemma}{Lemmata}
	\Crefname{cor}{Corollary}{Corollaries}
	\Crefname{defn}{Definition}{Definitions}
	\Crefname{rem}{Remark}{Remarks}
	\Crefname{exmp}{Example}{Examples}
	\Crefname{notation}{Notation}{Notations}
	\Crefname{ass}{Assumption}{Assumptions}

\numberwithin{table}{section}
\numberwithin{equation}{section}

\usepackage{array}
\numberwithin{subsection}{section}

\author{Michele Zordan}
\address{Fakult\"at f\"ur Mathematik, Universit\"at Bielefeld, D-33501 Bielefeld, Germany}
\email{zordan.michele@gmail.com}
\curraddr{KU Leuven, Departement Wiskunde, Celestijnenlaan 200B, B-3001 Leuven, Belgium.}

	\newcommand{\N}{\mathbb{N}}
	\newcommand{\Z}{\mathbb{Z}}
	\newcommand{\Q}{\mathbb{Q}}
	\newcommand{\C}{\mathbb{C}}

	\newcommand{\Mat}[2]{\mathop{\mathrm{Mat}_{#1}#2}}
	
	\DeclareMathOperator{\Hom}{Hom}
	\DeclareMathOperator{\coker}{coker}
	\newcommand{\im}{\mathop{\mathrm{im}}}
	\DeclareMathOperator{\linspan}{Span}
	\DeclareMathOperator{\tr}{tr}
	
	\newcommand{\rk}{\mathop{\mathrm{rk}}}
	
	\DeclareMathOperator{\id}{id}

	\newcommand{\field}{\mathbf{K}}
	\newcommand{\numfield}{k}
	\newcommand{\ringint}{\mathcal{O}}

		\newcommand{\completion}{\mathfrak{o}}
		\newcommand{\primeideal}{\mathfrak{p}}
		
		\newcommand{\qfield}{\mathbb{F}_q}
       			\newcommand{\level}{r}
        			\newcommand{\otherlevel}{t}

	\newcommand{\varY}{\mathbf{Y}}

		\DeclareMathOperator{\alggroup}{\mathbf{G}}

		\DeclareMathOperator{\lie}{Lie}

			\newcommand{\GL}[3]{\mathop{\mathrm{GL}_{#1}^{#3}{#2}}}
			\newcommand{\SL}[3]{\mathop{\mathrm{SL}_{#1}^{#3} #2}}
			
			\newcommand{\linaut}[1]{\mathrm{GL}{#1}}
			\newcommand{\spl}[2]{\mathop{\mathfrak{sl}_{#1}#2}}
			\newcommand{\gl}[2]{\mathop{\mathfrak{gl}_{#1}#2}}
			
			\DeclareMathOperator{\proj}{proj}
			\newcommand{\ad}[1]{\mathrm{ad}_{#1}}
			\newcommand{\Ad}[1]{\mathrm{Ad}_{#1}}
	\newcommand{\group}{G}
	\DeclareMathOperator{\padiclie}{L}

			\newcommand{\centr}[2]{\mathop{\mathrm{C}_{#1}\left(#2\right)}}
			\DeclareMathOperator{\stab}{Stab}
			
	\newcommand{\ring}{R}
		
			\DeclareMathOperator{\ch}{\mathrm{char}}
		
	\newcommand{\dimalg}{d}
	\newcommand{\liering}{\mathfrak{g}}
	\newcommand{\killing}{\kappa}
			\newcommand{\basis}{\mathcal{B}}
			\newcommand{\otherbasis}{\mathcal{H}}

	\newcommand{\zetafunc}[1]{\mathop{\zeta_{#1}(s)}}

		\newcommand{\halfdim}{h}
		\newcommand{\permiss}{m}
		\newcommand{\element}{a}
	\newcommand{\tuple}[1]{\mathrm{\mathbf{#1}}}

\newcommand{\size}{n}

		\newcommand{\form}{B}
	\newcommand{\reduction}[2]{\theta^{#1}_{#2}}

	\newcommand{\reduc}[2]{\theta_{#1,#2}}
	\newcommand{\Reduc}[2]{\Theta_{#1,#2}}
	\newcommand{\red}{\reduc{\level + 1}{\level}}
	\newcommand{\Red}{\Reduc{\level + 1}{\level}}
	\newcommand{\redp}[1]{\reduc{#1}{1}}
	\newcommand{\Redp}[1]{\Reduc{#1}{1}}
	\newcommand{\Redpimage}[1]{\Upsilon}

        \newcommand{\grp}[1]{G_{#1}}
        \newcommand{\alg}[1]{\mathfrak{g}_{#1}}
        \newcommand{\cimage}{\mathop{\widetilde{S}}}
        
        \newcommand{\kersigma}{\widetilde{N}}

	\newcommand{\alginf}{\liering_{}}
        \newcommand{\algfin}{\bar{\liering_{}}}

        \newcommand{\grpfin}{\bar{\group}}
        \newcommand{\grpinf}{\group}
        \newcommand{\Al}{\liering_{\level}}
        \newcommand{\Alone}{\liering_{\level+1}}
        \newcommand{\Grl}{\group_{\level}}
        \newcommand{\Grlone}{\group_{\level+1}}
        
	\newcommand{\pseries}[2]{\mathop{\mathcal{P}_{#1}#2}}
	\newcommand{\pcoeff}[2]{N_{#1,\mathrm{r}_{#1}}^{\completion}{#2}}
	\newcommand{\pcoeffCard}[2]{\mathrm{N}_{#1,\mathbf{r}_{#1}}^{\completion}{#2}}
		\newcommand{\Cmatrix}[2]{\mathop{\mathcal{R}_{#1}(#2)}}
		\newcommand{\cmatrix}{\mathcal{R}}
		\newcommand{\lastindex}{\ell}
		
	\newcommand{\cone}[1]{\mathop{W_{#1}^{\completion}}}
        
	\DeclareMathOperator{\Sh}{\mathrm{Sh}}
    	\newcommand{\sh}[2]{\mathop{\Sh_{#1}\mathopen{}\left({#2}\right)\mathclose{}}}
		\newcommand{\addspan}[1]{\mathop{\mathrm{As}(#1)}}
		\newcommand{\allshadows}[1]{\mathop{\mathfrak{Sh}(#1)}} 
		\newcommand{\shadowsequence}{\mathcal{I}}
		\newcommand{\coeff}[3]{\mathop{\mathcal{N}_{#1,#2}}}
		\newcommand{\shadowsequences}[1]{\mathop{\mathcal{D}_{#1}}}
		\newcommand{\isoshadow}{\mathbf{S}}
		\newcommand{\otherisoshadow}{\mathbf{T}}
		\newcommand{\ninshadow}[2]{\Lambda(#1,#2)}
        
	        \newcommand{\shadow}{S}
	        \newcommand{\othershadow}{T}
	
		\newcommand{\lieshadow}{\mathfrak{s}}

        		\newcommand{\orbit}{\mathcal{C}}
        		\newcommand{\liftorbit}{\widetilde{\mathcal{C}}}
		\newcommand{\halfdimshadow}[1]{\mathop{\delta(#1)}}
		\newcommand{\zen}[2]{\mathop{\mathrm{Triv}_{#1}(#2)}}
		\newcommand{\dimzen}[1]{\mathop{z_{#1}}}
		\newcommand{\poly}[2]{\mathop{\mathrm{f}_{#1}(#2)}}
		\newcommand{\npoly}[2]{\mathop{\mathrm{f}'_{#1}(#2)}}
		\newcommand{\minimal}[1]{m_{#1}}
		
		\newcommand{\exponents}[1]{\mathbf{r}_{#1}}
		\newcommand{\trivialshadow}{\grpfin}
		\newcommand{\SubregNilp}{\mathbf{J}}
		\newcommand{\SubregSem}{\mathbf{L}}
			\newcommand{\Reg}{\mathbf{R}}

		\newcommand{\zetasc}[2]{\zeta^\mathrm{sc}_{#1}#2}
		\newcommand{\nsc}[2]{a^\mathrm{sc}_{#1}{(#2)}}

\usepackage[disable]{todonotes}

	\usepackage[msc-links,initials]{amsrefs}

\title[Adjoint Orbits and Representation Zeta Functions]{Adjoint Orbits of Matrix Groups over Finite Quotients of Compact Discrete Valuation Rings and Representation Zeta Functions}

\begin{document}
\begin{abstract}
This paper gives methods to describe the adjoint orbits of $\alggroup(\completion_\level)$ on $\lie(\alggroup)(\completion_\level)$ where $\completion_\level=\completion/\primeideal^\level$ ($\level\in\N$) is a finite quotient of the completion $\completion$ of the ring of integers of a number field at a prime ideal $\primeideal$ and $\alggroup$ is a closed $\Z$-subgroup scheme of $\GL{\size}{}{}$ for an $\size\in\N$ such that the Lie ring $\lie(\alggroup)(\completion)$ is quadratic. The main result is a classification of the adjoint orbits in $\lie(\alggroup)(\completion_{\level+1})$ whose reduction $\bmod\,\primeideal^{\level}$ contains $\element\in\lie(\alggroup)(\completion_\level)$ in terms of the reduction $\bmod\,\primeideal$ of the stabilizer of $\element$ for the $\alggroup(\completion_\level)$-adjoint action. As an application, this result is then used to compute the representation zeta function of the principal congruence subgroups of $\SL{3}{(\completion)}{}$.
\end{abstract}
\maketitle
\thispagestyle{empty}
\section{Introduction}
\label{sec:intro}
\subsection{Main results} Let $\alggroup$ be a smooth closed $\Z$-subgroup scheme of $\GL{\size}{}{}$ for some $\size\in\N$. Let $\numfield$ be a number field with ring of integers $\ringint$. Let $\completion$ be the completion of $\ringint$ at a non-zero prime ideal  $\primeideal\vartriangleleft\ringint$ such that the map $\alggroup(\completion)\rightarrow\alggroup(\completion/\primeideal^\level)$ given by the reduction $\bmod\,\primeideal^\level$ is surjective for all $\level\in\N$. By Hensel's lemma this happens for all but finitely many prime ideals of $\ringint$ (see \cite[Chapter II, Proposition 4.1]{KNVS2011topics}). Let $\pi$ be a uniformizer for $\primeideal$ and identify the residue field $\completion/\primeideal$ with $\qfield$. For convenience of notation, in what follows we shall set $\completion_\level=\completion/\primeideal^\level$.

\begin{defn}
\label{shadow}
Let $\level\in\N$ and $\element\in\lie(\alggroup)(\completion_{\level})$. We define the (\emph{group}) \emph{shadow} 
\[\sh{\alggroup(\completion_{\level})}{a}\allowbreak\leq\alggroup(\qfield)\] 
of $a$ to be the reduction $\bmod$ $\primeideal$ of the group stabilizer of $\element$ for the adjoint action of $\alggroup(\completion_{\level})$ on $\lie(\alggroup)(\completion_{\level})$. Analogously, the \emph{Lie shadow} 
\[\sh{\lie(\alggroup)(\completion_{\level})}{\element}\leq\lie(\alggroup)(\qfield)\]
of $a$ is the reduction $\bmod$ $\primeideal$ of the centralizer of $\element$ in $\lie(\alggroup)(\completion_{\level})$.
\end{defn}
\begin{rem}
\Cref{shadow} borrows from \cite[Definition 2.2]{akov2}. The crucial difference here is that \cite[Definition 2.2]{akov2} also associates a conjugacy class of such shadows to each adjoint orbit in $\lie(\alggroup)(\completion_{\level})$. We shall work instead with individual elements.
\end{rem}
\begin{ass}
\label{ass:same_shadow}
For the rest of the section we fix $\level\in\N$ and $\element\in\lie(\alggroup)(\completion_{\level})$ having a lift to $\lie(\alggroup)(\completion_{\level+1})$ with the same shadow. We assume further that  $\lie(\alggroup)(\completion)$ is {\em quadratic}, i.e.\ it admits a non-degenerate ad-invariant symmetric form .
\end{ass}
The class of $\Z$-subgroup schemes such that $\lie(\alggroup)(\completion)$ is quadratic encompasses all semisimple algebraic groups defined over $\Z$ and such that the determinant of the Killing form on the associated Lie algebra is invertible in $\completion$. There are however examples that are not semisimple. An important one is $\alggroup = \GL{\size}{}{}$ with the form $\tr(XY)$ on $\lie(\alggroup)(\completion)$. Another comparatively easy example is the class-$2$ free nilpotent group on $3$ generators: it is the unipotent $\Z$-group scheme associated with the class-$2$ free nilpotent $\Z$-Lie lattice on $3$ generators $\mathfrak{n}_{3,2}$. The Lie lattice $\mathfrak{n}_{3,2}\otimes_\Z \completion$ is quadratic for almost all primes $\primeideal$ (see \cite[Theorem 6.1 (v)]{benconlal2016quadratic}).\par
The first main result concerns adjoint orbits in $\lie(\alggroup)(\completion_{\level})$.
\begin{main}
\label{thm:A}
Let $\level$, $a$ and $\alggroup$ be as in \cref{ass:same_shadow}. The set of $\alggroup(\completion_{\level+1})$-adjoint orbits in $\lie(\alggroup)(\completion_{\level+1})$ containing a lift of the element $\element$ is in one to one correspondence with the set of orbits for the co-adjoint action of $\sh{\alggroup(\completion_{\level})}{a}$ on $\Hom_{\qfield}(\sh{\lie(\alggroup)(\completion_{\level})}{\element},\qfield)$.
\end{main}
In case $\alggroup=\GL{n}{}{}$ and $\level=2$, Theorem A is \cite[Theorem 1]{jample2012normalforms}. Indeed, as proved in \cite[Lemma 6]{jample2012normalforms} for any $n\times n$ matrix over $\qfield$ there is an $n\times n$ matrix over $\completion_{2}$ with the same shadow lifting it. With the further hypothesis of the existence of a lift with the same shadow (cf.\ \cref{ass:same_shadow}), the proof of \cref{thm:A} generalizes the strategies adopted by Jambor and Plesken in \cite{jample2012normalforms}.\par
The second main result describes the shadow of a lift: 
\begin{main}
\label{cor:induction}
\label{thm:B}
Let $\level$, $a$ and $\alggroup$ be as in \cref{ass:same_shadow}. Let $x\in\lie(\alggroup)(\completion_{\level+1})$ be a lift of $\element\in\lie(\alggroup)(\completion_{\level})$, and let the orbit of $x$ for the action of $\alggroup(\completion_{\level+1})$ be represented by the orbit of 
\[c\in\Hom_{\qfield}(\sh{\lie(\alggroup)(\completion_{\level})}{\element},\qfield)\] 
in the one to one correspondence of \Cref{thm:A}. Then 
\[\sh{\alggroup(\completion_{\level+1})}{x} \cong \stab_{\sh{\alggroup(\completion_{\level})}{\element}}(c), \]
where $\stab_{\sh{\alggroup(\completion_{\level})}{\element}}(c)$ is the stabilizer of $c$ for the 
dual of the $\sh{\alggroup(\completion_{\level})}{\element}$-adjoint action on $\sh{\lie(\alggroup)(\completion_\level)}{\element}$.
\end{main}
The third main result is a quantitative statement about the number of lifts of a matrix. Let $\dimalg = \dim_{\qfield} \lie(\alggroup)(\qfield)$.
\begin{main}
\label{thm:D}
Let $\level$, $a$ and $\alggroup$ be as in \cref{ass:same_shadow}. Let $\shadow = \sh{\alggroup(\completion_\level)}{\element}$ and let $\othershadow$ be the shadow of a lift of $\element$ to $\lie(\alggroup)(\completion_{\level+1})$. Let $\lieshadow = \sh{\lie(\alggroup)(\completion_{\level})}{\element}$ and
\[\lambda=\lvert\lbrace c\in \Hom_{\qfield}(\lieshadow,\qfield)\mid \stab_{\shadow}{(c)}\cong\othershadow\rbrace\rvert,
\]
where $\stab_{\shadow}{(c)}$ is defined as in \cref{thm:B}. Then the number of lifts of $\element$ with shadow isomorphic to $\othershadow$ is equal to
\[
q^{\dimalg-\dim_{\qfield} \lieshadow}\, \lambda.
\]
\end{main}
The fourth main result is an application of the previous main results to representation zeta functions. Let $\group=\alggroup(\completion)$ have finite abelianization (FAb for short) i.e.\ $\lvert G/[G,G]\rvert < \infty$. By \cite[Proposition 2.1]{akov2013representation}, $\group$ is (representation) rigid i.e.\ the number $r_i(\group)$ of continuous complex $i$-dimensional irreducible representations is finite for each $i\in\N$, its \emph{representation zeta function} is the Dirichlet series
 \[
 \zetafunc{\group}=\sum_{i=1}^{\infty} r_i(G) i^{-s}\text{ ($s\in\C$).}
 \]
\Cref{thm:A,thm:B} are used to obtain the following result.
\begin{main}
\label{thm:C}
Let $\completion$ be a compact discrete valuation ring of characteristic $0$ whose residue field has cardinality $q>2$ and characteristic $p\neq 3$. Then for all $\permiss\in\N$ such that the $\permiss$-th principal congruence subgroup $\SL{3}{(\completion)}{\permiss}$  is potent and saturable (cf.\ \cite[Section 2.1]{akov2013representation}),
\[
\zetafunc{\SL{3}{(\completion)}{\permiss}}=q^{8\permiss}\frac{1+u(q)q^{-3-2s}+u(q^{-1})q^{-2-3s}+q^{-5-5s}}{(1-q^{1-2s})(1-q^{2-3s})}
\]
where $u(X)=X^3+X^2-X-1-X^{-1}$.
\end{main}
Here the $\permiss$-th principal congruence subgroup of $\alggroup(\completion)$ is the kernel of the reduction modulo $\primeideal^\permiss$, i.e.\
\[
\alggroup^\permiss(\completion) = \ker \left(\alggroup(\completion)\rightarrow\GL{n}{(\completion/\primeideal^\permiss)}{}\right).
\]

	\renewcommand{\level}{\ell}
This result already appeared as part of \cite[Theorem E]{akov2013representation} and was obtained again in \cite{akov2} by different methods. In \cite{akov2013representation} the representation zeta function is expressed as a Poincar\'e series, which is then computed with $\primeideal$-adic integration. In \cite{akov2} the authors give an expression of the representation zeta function in terms of certain shadow-similarity class zeta functions defined in \cite[Definition 5.14]{akov2}. Our approach is a hybrid of these previous two, i.e.\ we use shadows in order to compute the Poincar\'e series in \cite{akov2013representation}. There are three main ideas. The first one is that shadows may be related to kernels of a commutator matrix, this is shown in \cref{lem:Lie_induction} and exemplified in \cref{sec:subreg_nilp}. The second idea is new and is that the previous concept may be applied recursively to write a streamlined formula for the Poincar\'e series in \cite{akov2013representation}. This is the content of \cref{sec:poin_shadows} and culminates in \cref{eq:sum_products}. The last main idea is also new and is that \cref{eq:sum_products} may be further simplified discarding the differences among regular shadows (cf.\ \cref{sec:reg_convention}). See also \cref{rem:comparison} for a more detailed comparison between our methods and the ones in \cite{akov2}.
	\renewcommand{\completion}{A}
	\renewcommand{\primeideal}{I}
	\newcommand{\fooring}{A}
	
\subsection{Background and motivation}
In order to contextualize the main results of this paper and provide motivation for them, we now make a brief digression summarizing some known results on similarity classes.\par
When considering matrices over a field, the similarity classes are characterized by \emph{rational} (or \emph{Frobenius}) \emph{canonical forms} (see Dummit and Foote \cite[Section 12.2]{dumfoo2004algebra}). When the base ring is not a field, even over $\Z$ or its finite quotients, canonical forms are no longer available; nonetheless, over the years, many notable results have been proved.
In \cite{dav1968equations} Davis showed that, for a rational prime $p$ and $\level\in\N$, two matrices in $\Mat{n}{(\Z/p^\level \Z)}$, which are zeroes of a common polynomial whose reduction modulo $p$ has no repeated roots, are similar if and only if they are similar modulo $p$. 
In a similar flavour and generalizing a result of Suprunenko~\cite{sup1964conjugacy}, Pomfret showed that, over finite local rings, invertible matrices of order coprime to the residue field characteristic are similar if and only if their reductions modulo the maximal ideal are similar (see \cite{pom1973similarity}).\par
Another source of insights comes from the solution of the conjugacy problem for arithmetic groups achieved by Grunewald and  Segal. In \cite{gru1980conjugacy}  Grunewald gave a method to determine if two matrices in $\GL{n}{(\Q)}{}$ are conjugate by an invertible matrix over $\Z$. The same author and Segal, described in \cite{gruseg1979decision} a more general algorithm to decide whether two elements of an arithmetic group are conjugate.
For traceless $3\times3$ matrices over $\Z$, Appelgate and Onishi gave in \cite{apponi19823x3} an independent solution to the problem, giving a more effective algorithm to determine whether two matrices of $\SL{3}{(\Z)}{}$ are similar.\par
For $3\times 3$ matrices over $\Z/p^\level \Z$ ($\level\in\N$) \--- and slightly more generally over a finite quotient of a discrete valuation ring $\completion$ modulo a power of its maximal ideal $\primeideal$ \--- the first attempts of classifying the conjugacy classes date back at least to Nechaev~\cite{nec1983similarity}, where the similarity classes in $\Mat{3}{(\Z/p^2\Z)}$ are described. Pizarro in \cite{piz1983similarity} gave a complete classification for matrices over finite quotients of discrete valuation rings. More recently, in \cite{avnonnpravas2009similarity} Avni, Onn, Prasad and Vaserstein  have extended the classification of Nechaev classifying similarity classes of $3\times 3$ matrices over all finite quotients of $\completion$. This classification is explicit enough to allow them to enumerate the similarity classes in $\Mat{3}{(\completion/\primeideal^\level)}$ and the conjugacy classes of $\GL{3}{(\completion/\primeideal^\level)}{}$ for $\level\in\N$ (see \cite[Theorem 5.2]{avnonnpravas2009similarity}).\par
Even for $\level = 2$, the conjugacy problem for matrices in $\Mat{4n}{(\Z/p^2\Z)}$ contains, according to Nagorny{\u\i} \cite[Section 4]{nag1978complex}, the matrix pair similarity problem, which, according to Drozd \cite{dro1980tamewild}, is wild for general $n$. Nevertheless recent striking results have been obtained for similarity classes of matrices of arbitrary size over a local principal ideal ring of length $2$. Let $\ring$ be such a ring with residue field $\mathbf{F}$ of cardinality $t$. First, if $\ring'$ is another  local principal ideal ring of length $2$,  Singla  \cite{sin2010repr_gl} has shown that there is a canonical bijection between irreducible representations of $\GL{n}{(\ring)}{}$ and of $\GL{n}{(\ring')}{}$. In particular the number of conjugacy classes of these two groups is equal and only depends on the characteristic of the residue field. Second, Jambor and Plesken have proved that the similarity classes in $\Mat{n}{(\ring)}$ whose image over the residue field $\mathbf{F}$ is the similarity class of $\element\in\Mat{n}{(\mathbf{F})}$ are in one to one correspondence with the orbits of group centralizer $\centr{\GL{n}{(\mathbf{F})}{}}{\element}$ acting on the $\mathbf{F}$-linear dual of the commuting algebra $\centr{\Mat{n}{(\mathbf{F})}}{\element}$. More recently, Prasad, Singla and Spallone have formulated and proved an equivalent result phrased in terms of the $\mathrm{Ext}$ functor (see  \cite[Remark 1.1 and Theorem 2.8]{prapoospa2015similarity}). Using this theory, they describe the similarity classes in $\Mat{n}{(\ring)}$ for $n\leq 4$, together with their centralizers. This allows them to enumerate the similarity classes and the cardinalities of their centralizers as polynomials in $t$. In particular they show that the polynomials representing the number of similarity classes in $\Mat{n}{(\ring)}$ have non-negative integer coefficients.\par
\renewcommand{\alggroup}{\mathbf{\Gamma}}
\subsubsection{Zeta functions}
When $\ring = \completion/\primeideal^2$, the computations in \cite{prapoospa2015similarity} give the number of similarity classes of $\Mat{n}{(\ring/\primeideal^\level)}$ for $\level = 2$. If $\level$ is allowed to vary, natural questions on properties of the number of similarity classes of  $\Mat{n}{(\ring/\primeideal^\level)}$ as $\level$ tends to infinity arise. Slightly more generally, for an $\completion$-group scheme $\alggroup$, one studies the asymptotic behaviour of the number of $\alggroup(\completion/\primeideal^\level)$-adjoint orbits on the associated Lie lattice $\lie(\alggroup)(\completion/\primeideal^\level)$. Such questions may be addressed by means of the \emph{similarity class zeta function}
\[
\zetasc{\lie(\alggroup)(\completion)}{(s)}=\sum_{i\in\N} \nsc{i}{\lie(\alggroup)(\completion)} q^{-i s},
\]
where $\nsc{i}{\lie(\alggroup)(\completion)}$ denotes the number of $\alggroup(\completion/\primeideal^i)$-adjoint orbits in $\lie(\alggroup)(\completion/\primeideal^i)$ and $s$ is a complex variable. 
For odd residue field characteristic, Avni, Onn, Prasad and Vaserstein have computed $\zetasc{\gl{3}{(\completion)}}{}$ (cf.\  \cite[Theorem 5.2]{avnonnpravas2009similarity}) while Avni, Klopsch, Onn, and Voll in \cite[Theorem E]{akov2} computed the similarity class zeta function of $\mathfrak{gu}_{3}{(\completion)}$. Both similarity class zeta functions are rational in $q^{-s}$ and have abscissa of convergence $3$.\par
In a similar vein, du Sautoy \cite{sau2005conjugacy} proves that the zeta function counting conjugacy classes in congruence quotients of compact $p$-adic analytic groups is rational in $p^{-s}$. 
In particular this holds for $\GL{n}{(\Z_p)}{}$, establishing that there is a linear 
recurrence relation among the numbers of conjugacy classes of the groups 
$\GL{n}{(\Z/p^\level \Z)}{}$ ($\level\in\N$). 
More recently, Berman, 
Derakhshan, Onn and Paajanen have proved an analogous result for Chevalley groups over complete discrete valuation 
rings with sufficiently large residue field characteristic (see \cite[Theorem C]{berderonnpir2013uniform}). \par
	\renewcommand{\completion}{\mathfrak{o}}
	\renewcommand{\primeideal}{\mathfrak{p}}
	\renewcommand{\level}{r}
%
Another interesting application of classifying adjoint orbits in $\lie(\alggroup)(\completion)$ is computing  representation zeta functions. Classes of groups for which these have been studied so far comprise arithmetic groups and their principal congruence subgroups. For what concerns principal congruence subgroups of special linear groups, the Kirillov orbit method \--- when applicable \--- is a powerful linearization technique that relates irreducible representations and similarity classes. In \cite{akov2} Avni, Klopsch, Onn and Voll use the classification of adjoint orbits in $\gl{3}{(\completion)}$ and $\mathfrak{gu}_{3}(\completion)$ to compute the representation zeta function of principal congruence subgroups of $\SL{3}{(\completion)}{}$ and $\mathrm{SU}_{3}(\completion)$ in the same hypotheses of \cref{thm:C}.
\renewcommand{\alggroup}{\mathbf{G}}
\subsection{Organization of the paper}
We start off in \cref{sec:group_schemes} with a quick introduction to the vocabulary of group schemes over $\Z$, contextualizing this topic to the main purpose of the paper. We introduce a Lie theory for group schemes and the exponential map for closed subgroup schemes of $\GL{n}{}{}$. All results contained in this section are well known to the experts but difficult to find in the literature from a unique source; we therefore, for the sake of completeness, included them here.
\Cref{sec:adjoint_orbits} introduces our version of the similarity class invariant called the \emph{shadow}. We use it to generalize results of Jambor and Plesken (see \cite{jample2012normalforms}) and obtain \cref{thm:A,thm:B}, from which \cref{thm:D} is then deduced. The section ends with a refinement of \cref{thm:D} for special linear groups that is more suited to be used in the subsequent computations. 
\Cref{sec:appl}, finally, is concerned with applying the results in \cref{sec:adjoint_orbits} to the computation of representation zeta functions. 
	
\subsection{Notation}
We denote by $\N$ the set of the positive integers $\lbrace 1,2,\dots\rbrace$, while $\N_0=\lbrace 0,1,2,\dots\rbrace$ are the natural numbers. Analogously, for $n\in\N$ we set $[n]=\lbrace 1,\dots,n\rbrace$ and $[n]_0=\lbrace 0,\dots,n\rbrace$. In this work, $p$ is a rational prime. The field of $p$-adic numbers is denoted by $\Q_p$ and the ring of $p$-adic integers by $\Z_p$. \par
The group of units of a ring $\ring$ is $\ring^*$. We introduce a similar notation for non-trivial $\completion$-modules as follows. Given such a module $M$, we write $M^*=M\smallsetminus \primeideal M$. For the trivial $\completion$-module we set $\lbrace 0\rbrace^*=\lbrace 0\rbrace$.\par
If $\ring$ is a ring we write $R\llbracket T\rrbracket$ for the ring of formal power series in $T$. For $m\in\N$ and $f\in R\llbracket T\rrbracket$, $f\,\bmod T^m$ denotes the class of $f$ in the quotient ring $R\llbracket T\rrbracket/T^m$.
\subsection{Acknowledgements }
I am indebted to Christopher Voll and Benjamin Martin for their precious advice. I also wish to thank Tobias Rossmann, Giovanna Carnovale, Andrea Lucchini, Uri Onn and Alexander Stasinski for the interesting conversations and insightful comments on this work. I also wish to thank the referee for their useful comments on the first version of this paper.\par
This work is part of my PhD thesis. I acknowledge financial support from the School of Mathematics of the University of Southampton, the Faculty of Mathematics of the University of Bielefeld and CRC 701. I am currently supported by the Research Project G079218N of the Research Foundation - Flanders (FWO).\par

\section{Group Schemes}
\label{sec:group_schemes}
\subsection{Preliminaries on group schemes}
\label{sec:group_schemes/prel}
An (\emph{affine}) \emph{group scheme} $\alggroup$ over $\Z$ (or $\Z$-group scheme) is a $\Z$-group functor that is representable when considered as a functor from $\mathbf{Rng}$ to $\mathbf{Set}$. A $\Z$-subgroup scheme is a subscheme of a $\Z$-group functor that is also a group scheme in its own right. What follows is a summary of some basic concepts in the theory of groups schemes. We refer to \cite{wat1979introduction} for an introduction to group schemes and to \cite{demgab1980alggroups} for a more advanced treatment. 
\subsubsection{The Lie algebra of a group functor}
Let  $\ring$ be a ring, if $\ring[T]$ is the algebra of polynomials in $T$ with coefficients in $\ring$, we write $\varepsilon$ for the class of $T\,\bmod\,T^2$ and $\ring(\varepsilon)$ for the quotient algebra $\ring[T]/T^2$. We have a decomposition $\ring(\varepsilon) = \ring\oplus \varepsilon \ring$ and homomorphisms $i:\ring\rightarrow\ring(\varepsilon)$, $\mathrm{\proj}:\ring(\varepsilon)\rightarrow\ring$ defined by $i(1) = 1$ and $\mathrm{proj}(1) = 1$,  $\mathrm{proj}(\varepsilon) = 0$, such that $\mathrm{proj}\circ i = \id_\ring$.\par
Let $\alggroup$ be a group scheme over $\Z$. The homomorphisms $i$ and $\mathrm{proj}$ define homomorphisms $\alggroup(i):\alggroup(\ring)\rightarrow \alggroup(\ring(\varepsilon))$ and $\alggroup(\mathrm{proj}): \alggroup(\ring(\varepsilon))\rightarrow \alggroup(\ring)$. The $\Z$-group functor $\lie(\alggroup)$ is defined by
\[
\lie(\alggroup)(\ring) = \ker (\alggroup(\mathrm{proj})).
\]
By \cite[II, \S 4, 4.1]{demgab1980alggroups} $\lie(\alggroup)(\ring)$ has the structure of an $\ring$-lattice (i.e.\ a free $\ring$-module of finite rank). When no risk of confusion exists, by abuse of notation, $\alggroup(i)$ and $\alggroup(\mathrm{proj})$ will also be denoted by $i$ and $\mathrm{proj}$. 
\subsubsection{The linear group}
We now introduce a very important example of $\Z$-group scheme, namely the (general) linear group. If $V$ is a $\Z$-module (i.e.\ an abelian group) and $\ring$ is a ring, $\mathcal{L}(V\otimes_{\Z} \ring)$ denotes the monoid of all endomorphisms of the $\ring$-module $V\otimes_{\Z} \ring$. We define a $\Z$-monoid functor $\mathrm{End}(V)$ by setting
\[
\mathrm{End}(V)(S) = \mathcal{L}(V\otimes_{\Z} S)\text{ ($S\in\mathbf{Rng}$)}.
\]
The \emph{linear group} of $V$, denoted by $\GL{}{(V)}{}$, is the largest subgroup functor of $\mathrm{End}(V)$. The discussion in \cite[II,\S, 2.4]{demgab1980alggroups} shows that $\mathrm{End}(V)$ is an affine scheme over $\Z$ when $V$ is a free finitely generated abelian group. If $V = \Z^\size$, we write $\GL{\size}{}{} = \GL{}{(V)}{}$.
\subsubsection{The adjoint action}
The group $\alggroup(\ring)$ acts on $\lie(\alggroup)(\ring)$ in the following way: let $g$ be an element of $\alggroup(\ring)$ and $x\in \lie(\alggroup)(\ring)$, we set
\[
\Ad{g}(x) = i(g)\,x\,i(g)^{-1}.
\]
Writing $\linaut{(\lie(\alggroup))}$ for $\GL{}{(\lie(\alggroup)(\Z))}{}$ we may define a homomorphism
\[
\Ad{}: \alggroup\rightarrow \linaut{(\lie(\alggroup))}
\]
which is called the adjoint action of $\alggroup$. This in turn defines a homomorphism
\[
\ad{} : \lie(\alggroup)\rightarrow\lie(\linaut{(\lie(\alggroup))})
\]
by means of which one defines $[x,y] = \ad{}(x)(y)$ for all $x,y \in\lie(\alggroup)(\ring)$. This gives $\lie(\alggroup)(\ring)$ the structure of an $\ring$-Lie lattice (i.e.\ an $\ring$-Lie ring that is also a free $\ring$-module of finite rank). For convenience of notation, we shall write
\[
 \ad{x} = \ad{}(x): \lie(\alggroup)(\ring)\rightarrow\lie(\alggroup)(\ring).
\]
\subsection{Exponential map}
\label{sec:group_schemes/exp_map}
The goal of this section is to introduce the exponential map on $\lie(\alggroup)(\completion)$. We need some notational conventions first.\par
It is customary to write the group law of $\lie(\alggroup)$ additively; we inherit the following notation from \cite{demgab1980alggroups}. If $S$ is an $\ring$-algebra and $\alpha$ is an element of $S$ of vanishing square, then there is a unique $\ring$-algebra homomorphism $\ring(\varepsilon)\rightarrow S$ sending $\varepsilon$ onto $\alpha$. The image of $x \in \lie(G)(R)$ under the composite homomorphism
\[
\lie(G)(R)\rightarrow \alggroup(\ring(\varepsilon))\rightarrow \alggroup(S)
\]
will be written $e^{\alpha x}$. Thus in $\alggroup(S)$ we have $e^{\alpha(x+y)} = e^{\alpha x} e^{\alpha y}$ for $x,y\in\lie(\alggroup)(\ring)$.\par
The following proposition introduces the exponential map in characteristic $0$ and is inspired by \cite[ II, \S 6, 3.1]{demgab1980alggroups}. We borrow from there the following convention: given a linearly topologized and complete $\ring$-algebra $S$, and a topologically nilpotent element $t$ of $S$, we write $f(t)$ for the element of $\alggroup(S)$ which is the image of $f(T)\in\alggroup(\ring\llbracket T\rrbracket)$ under the continuous morphism of $\ring\llbracket T\rrbracket$ into $S$ sending $T$ onto $t$. Therefore for instance, we shall have $f(\varepsilon)$ in $\alggroup(\ring(\varepsilon))$ and  $f(T + T')$ in $\alggroup(\ring\llbracket T, T'\rrbracket)$.
\begin{prop}
\label{prop:exp}
Let $\ring$ be a ring with $\ch\ring = 0$ and let $\alggroup$ be an affine group scheme over $\Z$. Then for each $x\in\lie(\alggroup)(\ring)$ there is  a unique element $\exp(Tx)$ of $\alggroup(\ring\llbracket T\rrbracket)$ such that
\begin{enumerate}
\item\label{exp_a} $\exp(\varepsilon x) = e^{\varepsilon x}$ in $\alggroup(\ring(\varepsilon))$,
\item\label{exp_b} $\exp((T + T') x) = \exp(T x) \exp(T' x)$ in $\alggroup(\ring\llbracket T, T'\rrbracket)$.
\end{enumerate}
\begin{proof}
An analogue of this result is proved in \cite[ II, \S 6, 3.1]{demgab1980alggroups} when  $\alggroup$ is a (not necessarily affine) group scheme over a field $\field$ of characteristic $0$ and $\ring$ is a $\field$-algebra. The argument there only uses that the base ring is a field to deduce that the group scheme is separated. Since  $\alggroup$ is affine here, it is separated as a group scheme over $\Z$. The rest of the proof goes through mutatis mutandis as in  loc.\ cit.\
\end{proof}
\end{prop}
\subsubsection{Closed subgroups of the linear group}
\label{sec:group_schemes/lin}
We shall now focus on a particular type of group schemes: closed subgroup schemes of the linear group. From this point onwards $\alggroup$ denotes a smooth closed $\Z$-subgroup scheme of $\GL{\size}{}{}$ ($\size\in\N$). 
By \cite[II, \S4, 4.12]{demgab1980alggroups} and references therein, $\lie(\alggroup)(\ring)$ may be identified with the $\ring$-Lie sublattice of all $x\in\mathrm{End}(\Z^\size)(\ring)$ such that $\id + \varepsilon x \in \alggroup (\ring(\varepsilon))$. With this identification, the adjoint representation is given by 
\[
\Ad{g}(x) = g\circ x\circ g^{-1} \text{ ($g\in\alggroup(\ring)$, $x\in\lie(\alggroup)(\ring)$),}
\]
and Lie bracket is the usual commutator of two matrices.
\begin{rem}
\label{rem:exp_gl}
For each $x\in\lie(\alggroup)(\completion)$,
\[
\exp(Tx) = \sum_{i\geq 0}\frac{T^i x^i}{ i ! }.
\]
On the right-hand side, $x^i$ ($i\in\N$) denotes the $i$-fold matrix multiplication of $x$ with itself. Notice, moreover, that \cite[Lemma 6.20]{dixsauman1999analytic} ensures that it makes sense to define the formal power series in $T$ on the right-hand side of the equality above.
\begin{proof}
Same as \cite[II, \S 6, 3.3]{demgab1980alggroups} with the fact that $\alggroup$ is closed and smooth.
\end{proof}
\end{rem}
\begin{rem}
\label{rem:exp_on_quotient}
Let $\level\in\N$. The exponential map $\exp:\lie(\alggroup)(\completion)\rightarrow\alggroup(\completion\llbracket T\rrbracket)$ induces an exponential map
\[
\exp_\level:\lie(\alggroup)(\completion_\level)\rightarrow\alggroup(\completion_\level\llbracket T\rrbracket).
\]
Namely, for each $x\in\lie(\alggroup)(\completion_\level)$ there is a unique element $\exp_\level(T x)$ such that
\begin{enumerate}
\item $\exp_\level(\varepsilon x) = e^{\varepsilon x}$ in $\alggroup(\completion_\level(\varepsilon))$,
\item $\exp_\level((T + T') x) = \exp_\level(T x) \exp_\level(T' x)$ in $\alggroup(\completion_\level\llbracket T, T'\rrbracket)$.
\end{enumerate}
\end{rem}
For convenience of notation, when there is no risk of confusion, we shall denote $\exp_\level$ with $\exp$ as well. Following the same arguments contained in \cite[II, \S 6, 3.4]{demgab1980alggroups}, the uniqueness statement in \cref{rem:exp_on_quotient} implies the following corollary.
\begin{cor}
\label{cor:exp_and_morphisms}
Let $x\in\lie(\alggroup)(\completion_\level)$. %
Then in $\GL{}{(\lie(\alggroup)(\Z))}{}(\completion_\level\llbracket T\rrbracket)$ we have
\[
\Ad{\exp(Tx)} = \sum_{i\geq 0} \frac{T^i \ad{x}^i}{ i !}.
\]
\end{cor}
When $x\in\pi\lie(\alggroup)(\completion_\level)$, then $\exp(Tx)\in \alggroup (\completion_\level[T])$ and we define $\exp(x)$ as the image of $\exp(Tx)$ under the homomorphism $\completion_\level[T]\rightarrow \completion_\level$ sending $T$ to $1$. In practice we may, then, replace $T$ by $1$ in \cref{cor:exp_and_morphisms} obtaining
\begin{equation}
\label{eq:Ad_ad}
\Ad{\exp(x)} = \sum_{i\geq 0} \frac{\ad{x}^i}{ i! }.
\end{equation}
\begin{prop}
\label{prop:exp_determines_lie_alg}
If $x\in\pi\gl{\size}{(\completion_\level)}$ then $x\in\pi\lie(\alggroup)(\completion_\level)$ if and only if $
\exp(x)\in\alggroup(\completion_\level)$.
\begin{proof}
We start by observing that if $x\in\pi\lie(\alggroup)(\completion_\level)$, then $\exp(Tx)\in\alggroup(\completion_\level[T])$, so we may replace $T$ with $1$ obtaining $\exp(x)\in\alggroup(\completion_\level)$. 
Moreover, by taking the Weil restriction of $\alggroup$ (cf.\ \cite[Section A.5]{congabpra2015pseudoreductive}) we may assume, within the proof of this statement, that $\completion = \Z_p$.\par
Let $\otherlevel$ be the valuation of the entry of $x$ with the lowest valuation. We observe that, as $\alggroup(\completion_\level)$ is a group and $\exp(x)\in\alggroup(\completion_\level)$, $\exp(x)^{p^{\level-\otherlevel-1}}\in\alggroup(\completion_\level)$.\par
Now, $\completion_\level \cong \Z/p^\level\Z$ and therefore, for each $\bar{a}\in \completion_\level$  and $a\in\Z$ such that $a\equiv \bar{a}\mod p^\level$, $
\exp(x)^{a} = \exp(\bar{a} x)$. 
This implies that $\exp(p^{\level-\otherlevel-1} x) = \id + p^{\level-\otherlevel-1} x \in \alggroup(\completion_\level)$, or equivalently that $p^{\level-\otherlevel-1} x\in\lie(\alggroup)(\completion_\level)$. The latter is equivalent to $x\in\lie(\alggroup)(\completion_\level)$ and we conclude.
\end{proof}
\end{prop}
\section{Adjoint Orbits in Lie Lattices}
\label{sec:adjoint_orbits}
Let $\alggroup$ be a closed smooth $\Z$-subgroup scheme of $\GL{\size}{}{}$ as before. We set $\grpinf=\alggroup(\completion)$ and $\alginf=\lie(\alggroup)(\completion)$; analogously, for all $\level\in\N$, $\grp{\level}=\alggroup(\completion_{\level})$ and $\alg{\level}=\lie(\alggroup)(\completion_{\level})$. For convenience of notation we write $\grpfin = \grp{1} = \alggroup(\qfield)$ and $\algfin = \alg{1} = \lie(\alggroup)(\qfield)$.\par
We set the following notation: for $\level,\otherlevel\in\N$ with $\level<\otherlevel$, we define 
\begin{align*}
\reduction{}{\level}:			& \alginf \rightarrow \alg{\level}			&& \level\in\N\\
\reduc{\level}{\otherlevel}:		&\alg{\level}\rightarrow \alg{\otherlevel} 	&& \level > \otherlevel
\end{align*}
to be the maps defined by reducing modulo $\primeideal^\level$ and $\primeideal^\otherlevel$ respectively. In a similar fashion $\Reduc{\otherlevel}{\level}$ denotes the reduction modulo $\primeideal^\level$ on $\grp{\otherlevel}$.
If  $a\in \alg{\otherlevel}$ for some $\otherlevel\in\N$, we say that $b\in\reduc{\level}{\otherlevel}^{-1}(a)$ is a \emph{lift} of 
$a$ to $\alg{\level}$. 
	\subsection{Shadows}
	\label{sec:shadows}
	Fix $\level\in\N$ and $\element\in\Al$. We denote the group and the Lie centralizer of $\element$ with
\begin{align*}
\centr{\Grl}{\element}&=\lbrace g\in\Grl\mid g\element g^{-1}=\element\rbrace\\
\centr{\Al}{\element}&=\lbrace x\in\Al\mid [x,\element]= 0 \rbrace
\end{align*}
respectively. 
%
The $\centr{\Grl}{\element}$-conjugation on $\centr{\Al}{\element}$ induces a $\sh{\Grl}{\element}$-action by conjugation on $\sh{\Al}{\element}$.
If $b\in \red^{-1}(\element)$ and $\liftorbit$ denotes its $\Grlone$-orbit, then $\liftorbit\cap \red^{-1}(\element)$ is an orbit of the action of 
\[\cimage=\Red^{-1}(\centr{\Grl}{\element})\]
on $\red^{-1}(\element)$. To see this, let $g\in \Grlone$ be such that $g.b\in \red^{-1}(\element)$, it follows that $\red(g.b)=\Red(g).\element=\element$, which implies $g\in\cimage$. 
	\subsection{Action of the first principal congruence subgroup}
	\label{sec:action_kernel}
	Let $\Redpimage{\level + 1}$ be the restriction of $\Redp{\level + 1}$ to $\cimage$. Following the approach of \cite{jample2012normalforms}, we proceed in two stages: first we consider the orbits for the action of the normal subgroup $\kersigma=\ker{\Redpimage{\level+1}}\trianglelefteq \cimage$ and then we act on them with the factor group $\cimage/\kersigma = \sh{\Grl}{\element}$. The following analogue of \cite[Lemma 5]{jample2012normalforms} describes the $\kersigma$-orbits in $\red^{-1}(\element)$.
\begin{lem}
\label{lem:action_kernel}
Let $b\in\red^{-1}(\element)$. Then there is a  one to one correspondence between $\pi^\level \coker \pi\ad{b}$ and the $\kersigma$-orbits in $\red^{-1}(\element)$.
\begin{proof}
We describe the $\kersigma$-conjugation in $\red^{-1}(\element)$ in terms of $b$ and the image of $\ad{b}$. We start by writing elements of $\red^{-1}(\element)$ in terms of $b$. Indeed, the latter is a preimage of $\element$ for the map $\red$. As any other preimage of $\element$ differs from $b$ by an element that is $0$ modulo $\primeideal^\level$, it follows that 
\begin{equation*}
	\red^{-1}(\element)=\lbrace b+\pi^{\level} z\mid z\in\Alone\rbrace.
\end{equation*}
By \cref{prop:exp_determines_lie_alg} an element of $\kersigma$ is of the form $\exp(\pi y)$ for some $\pi y\in \pi \Alone$. Now we are able to explicitly describe the $\kersigma$-conjugation in $\red^{-1}(\element)$. Fix $z\in\Alone$ and  $y\in \pi \Alone$. Let also $x = b+\pi^{\level} z\in\red^{-1}(\element)$ and $g = \exp(\pi y)$. Then
 \allowdisplaybreaks
\begin{align*}
g x g^{-1}  & = \Ad{\exp(\pi y)} (x)\\
		& = \sum_{i\geq 0}\frac{\ad{\pi y}^i(x)}{ i! }													&\text{ by \cref{eq:Ad_ad}}\\
		& = b + \pi^\level z + \sum_{i\geq 1}\left(\frac{\ad{\pi y}^i(b)}{ i! } + \frac{\ad{\pi y}^i(\pi^\level z)}{ i! }\right)\\
		& = b + \pi^\level z + \sum_{i\geq 1}\frac{\ad{\pi y}^i(b)}{ i! } 									& \text{ ($[\pi y,\pi^\level z]\equiv 0 \mod \primeideal^{\level+1}$).}
\end{align*}
We now show that $\ad{\pi y}^i(b) \equiv 0\,\bmod\,\primeideal^{\level+1}$, for $i\geq 2$.  To prove this, recall that $\exp(\pi y)\in \cimage = \red^{-1}(\centr{\Grl}{\element})$, therefore
\begin{equation*}
\Ad{\exp(\pi y)} (b)\equiv \sum_{i\geq 0}  \frac{\ad{\pi y}^i(b)}{ i! }   \equiv b \mod\primeideal ^\level.
\end{equation*}
Since this happens if and only if $\ad{\pi y}(b)\equiv 0\,\bmod\,\primeideal^\level$, the claim follows. This implies that $g x g^{-1} = b+\pi^{\level}z+ \pi [y,b]$, and shows that two lifts of $\element$, say $b+\pi^{\level} z$ and $b+\pi^{\level} z'$ ($z,z'\in\Alone$), are $\kersigma$-conjugate if and only if there exists $v\in\Alone$ such that $\pi^{\level} (z-z') = \pi [b,v]$; which suffices to conclude.
\end{proof}
\end{lem}
\begin{rem}
\label{rem:change_basis}
Let $s,t\in\N$ with $s<t$ and let $M$ be a left $\completion_t$-module, then $\pi^{t-s} M$ may be viewed as an $\completion_s$-module. Indeed, for all $\alpha\in\completion_{s}$ and $x\in M$, we define $\alpha x = \widehat{\alpha} x$, where $\widehat{\alpha}\in\completion_{t}$ is a lift of $\alpha$. This definition is unambiguous for the left multiplication of elements of $\pi^{t-s} M$ by elements of $\primeideal^{s}$ results in $0$.
\end{rem} 
In view of the last remark, we can formulate and prove the following lemma.
\begin{lem}
\label{lem:reduction}
Let $b$ be as in \cref{lem:action_kernel}. Then $\pi^{\level} \coker\pi\ad{b}\cong \pi^{\level -1}\coker\ad{a}$ as $\qfield$-vector spaces.
\begin{proof}
Let
\[
\Phi_\level:\xymatrix@R=3pt{\pi\gl{\size}{(\completion_{\level + 1})}\ar[r]&\gl{\size}{(\completion_{\level})}\\
\pi u\ar@{|->}[r]&u',}
\]
where $u'$ is the reduction modulo $\primeideal^\level$ of $u$, and let $\varphi_\level$ be the restriction of $\Phi_\level$ to $\pi\Alone$. The map $\varphi_\level:\pi\Alone\rightarrow\Al$ defines an isomorphism of $\completion_\level$-modules. As $\varphi_\level (\im \pi\ad{b}) = \im \ad{a}$, we have that $\varphi_\level$ induces an isomorphism $\bar{\varphi}_\level$ of $\qfield$-vector spaces between $\pi^{\level} \coker\pi\ad{b}$ and $\pi^{\level -1}\coker\ad{a}$.
\end{proof}
\end{lem}
\begin{notation}
\label{not:phi}
For further usage, we fix the the name $\varphi_\level$ for the restriction to $\pi\Alone$ of the map $\Phi_\level$ defined in the proof of \cref{lem:reduction} and we denote with $\bar{\varphi}_\level$ the $\qfield$-linear isomorphism between $\pi^{\level} \coker\pi\ad{b}$ and $\pi^{\level -1}\coker\ad{a}$ induced by $\varphi_\level$ as explained in the proof of \cref{lem:reduction}.
\end{notation}
\cref{lem:reduction} allows us to substitute $\pi^{\level} \coker\pi\ad{b}$ with $\pi^{\level-1} \coker \ad{a}$ on which $\sh{\Grl}{\element}$ acts with the action induced by the bijection  $\bar{\varphi}_\level$.
	\subsection{Action of the factor group}
	\label{sec:factor_group}
	We shall now investigate the action of the factor group $\cimage/\kersigma=\sh{\Grl}{\element}$ on the set of orbits for the $\kersigma$-action on $\red^{-1}(\element)$; i.e.\ we shall describe the action of  $\sh{\Grl}{\element}$ on $\pi^{\level -1}\coker\ad{a}$.
\begin{defn}
\label{def:conj}
The centralizer $\centr{\Grl}{\element}$ acts naturally by conjugation on $\pi^{\level-1} A$. Since $\exp(\pi\Al)\cap\centr{\Grl}{\element}$ is in its kernel, this action induces a $\sh{\Grl}{\element}$-action on $\pi^{\level-1} \Al$; namely an element $c\in\sh{\Grl}{\element}$ acts on $\pi^{\level-1} \Al$ conjugating by any of its lifts to $\centr{\Grl}{\element}$. We call this the $\sh{\Grl}{\element}$-\emph{conjugation by lifts}. Explicitly, if $\bar{c}\in\sh{\Grl}{\element}$ and  $c\in\centr{\Grl}{\element}$ is a lift of $\bar{c}$, for all $x\in\pi^{\level-1}\Al$, we write
\[
\bar{c}.x = c x c^{-1}.
\]
\end{defn}
The $\sh{\Grl}{\element}$-conjugation by lifts on $\pi^{\level-1}\Al$ induces an action on $\pi^{\level-1}\coker\ad{\element}$. Indeed, let $y\in\pi^{\level-1}\Al$, $\bar{c}\in\sh{\Grl}{\element}$ and let $c\in\centr{\Grl}{\element}$ be a lift of $\bar{c}$. As $c$ commutes with $\element$ 
\[
\bar{c}.[\element,y] = c (ay-ya) c^{-1} = \element cyc^{-1} - cyc^{-1} \element = [\element, \bar{c}.y].
\]
This implies that, denoting by $\Gamma_{\level,\bar{c}}$ the linear automorphism of $\pi^{\level-1}\Al$ defined by $x\mapsto \bar{c}.x$, and by $\rho_\level$ the projection of $\pi^{\level-1}\Al$ onto 
\[
\pi^{\level-1}(\Al/\im \ad{\element} ) =   \pi^{\level-1}\coker\ad{\element},
\]
there is a uniquely well defined $\qfield$-linear endomorphism $\bar{\Gamma}_{\level,\bar{c}}$ of $\pi^{\level-1}\coker\ad{\element}$ that makes the following diagram commute
\[
\xymatrix@R=30pt
{\pi^{\level-1}\Al				  \ar[r]^{\Gamma_{\level,\bar{c}}} \ar[d]_{\rho_\level}& \pi^{\level-1}\Al\ar[d]^{\rho_\level}\\
 \pi^{\level-1}\coker\ad{\element} \ar@{->}[r]^{\bar{\Gamma}_{\level,\bar{c}}}&         \pi^{\level-1}\coker\ad{\element}.}
\]
The rule $\bar{c}\mapsto \bar{\Gamma}_{\level,\bar{c}}$  defines a $\sh{\Grl}{\element}$-action on $\pi^{\level-1}\coker\ad{\element}$.\par
We shall now show that the $\sh{\Grl}{\element}$-action on $\pi^{\level-1}\coker\ad{\element}$ induced by $\bar{\varphi}_\level$ and resulting from the action of $\cimage/\kersigma = \sh{\Grl}{\element}$ on the set of orbits of the $\kersigma$-conjugation in $\red^{-1}(\element)$ is indeed the $\sh{\Grl}{\element}$-action on $\pi^{\level-1}\coker\ad{\element}$ described above. Analogously to the approach of \cite[Section 2.2]{jample2012normalforms}, the key to do this is finding a lift $b$ of $a$ with the same shadow. What we mean is made precise in the following definitions:
\begin{defn}
\label{def:same_shadow}
Let $\level\in\N$. We say that $b\in\Alone$ is \emph{shadow-preserving lift} of $a$ when $\red(b)=a$ and $\sh{\Grlone}{b}=\sh{\Grl}{a}$. We say that $\alg{}$ is \emph{shadow-preserving} if, for every $\level\in\N$, every $x\in\alg{\level}$ admits a shadow-preserving lift.
\end{defn}
\begin{exmp}
\label{ex:hered}
By  \cite[Lemma 6.4]{akov2}, the Lie lattice $\spl{3}{(\completion)}$ is shadow-preserving. 
\end{exmp}
The next lemma achieves what discussed above.
\begin{lem}
\label{cor:linear}
Assume that the element $\element$ admits a shadow-preserving lift. Then the action of $\sh{\Grl}{a}$ on $\pi^{\level -1}\coker\ad{a}$ induced by $\bar{\varphi_\level}$ is the linear action induced by the $\sh{\Grl}{a}$-conjugation by lifts.
\begin{proof}
Let $b\in\Alone$ be a shadow-preserving lift of $\element$. Analogously to \cref{def:conj} the group $\sh{\Grlone}{b}$ acts  on $\pi^\level\Alone$ by conjugation by lifts; as $b$ is shadow-preserving, this becomes an action of $\sh{\Grl}{\element}$ and it induces a $\sh{\Grl}{\element}$-action on $\pi^{\level}\coker\pi\ad{b}$ in the same way as $\sh{\Grl}{\element}$ induces an action on $\pi^{\level-1}\coker\ad{\element}$. These two actions commute with $\bar{\varphi}_\level$; in other words, if $\bar{c}\in\sh{\Grl}{\element}$,  the action by $\bar{c}$ on $\pi^\level\Alone$ defines $\qfield$-linear automorphisms $\Gamma_{\level+1,\bar{c}}$ of $\pi^{\level}\Alone$ and all cells in the following diagram commute
	\begin{equation}
	\label{eq:diag_induced_actions}
	\xymatrix@R=30pt
	{\pi^{\level}\coker\pi\ad{b} 		 	\ar[r]^{\bar{\Gamma}_{\level+1,\bar{c}}}&         								\pi^{\level}\coker\pi\ad{b}\\
	\pi^{\level}\Alone	        		 	\ar[r]^{\Gamma_{\level+1,\bar{c}}} \ar[u]^{\rho_{\level+1}}\ar[d]_{\varphi_\level}& 	\pi^{\level}\Alone\ar[u]_{\rho_{\level+1}}\ar[d]^{\varphi_\level}\\
	\pi^{\level-1}\Al		        		 	\ar[r]^{\Gamma_{\level,\bar{c}}} \ar[d]_{\rho_\level}& 							\pi^{\level-1}\Al\ar[d]^{\rho_\level}\\
	\pi^{\level-1}\coker\ad{\element} 	\ar[r]^{\bar{\Gamma}_{\level,\bar{c}}}&         								\pi^{\level-1}\coker\ad{\element}.}
	\end{equation}
Here $\rho_{\level+1}$ is the projection of $\pi^\level\Alone$ onto $\pi^{\level}\coker\pi\ad{b}$ and $\bar{\Gamma}_{\level+1,\bar{c}}$ is the $\qfield$-linear automorphism induced by $\Gamma_{\level+1,\bar{c}}$.\par
It follows that it suffices to prove that the $\sh{\Grl}{\element}$-conjugation by lifts on $\pi^\level\Alone$ induces the $\sh{\Grl}{\element}$-action on $\pi^{\level}\coker\pi\ad{b}$ obtained by letting $\cimage/\kersigma = \sh{\Grl}{\element}$ act on the set of orbits of the $\kersigma$-conjugation in $\red^{-1}(\element)$. Let $c\in\sh{\Grl}{a}$. Since $b$ has the same shadow as $a$, we can choose $\tilde{c}\in \centr{\Grlone}{b}$ lifting $c$. In order to see how $\tilde{c}$ acts on $\pi^{\level}\coker\pi\ad{b}$, first we see how it acts on an arbitrary lift of $a$:
\[
\tilde{c}(b+\pi^{\level} x)\tilde{c}^{-1}=b+\pi^{\level} \tilde{c} x\tilde{c}^{-1}.
\]
This last equation and \cref{lem:action_kernel} imply that the orbit of $\tilde{c}(b+\pi^{\level} x)\tilde{c}^{-1}$ corresponds to the class of $\pi^{\level} \tilde{c} x\tilde{c}^{-1}$ in $\pi^{\level}\coker\pi\ad{b}$. By \cref{eq:diag_induced_actions}, this allows us to conclude.
\end{proof}
\end{lem}
	\subsection{Intrinsic description of the orbits}
	\label{sec:intrinsic}
	 So far we have established a one to one correspondence between the $\Grlone$-orbits in $\Alone$ intersecting $\red^{-1}(\element)$ non-trivially and $\sh{\Grl}{a}$-conjugacy orbits in $\pi^{\level-1}\coker \ad{\element}$. Now we  replace $\pi^{\level-1}\coker \ad{\element}$ with the more intrinsic dual of the Lie shadow.
\begin{notation}
Let $\otherlevel\in\N$. Given an $\completion_\otherlevel$-module $M$ we write $M^\sharp$ for its dual, i.e. $M^{\sharp}=\Hom_{\completion_\otherlevel}(M,\completion_\otherlevel)$. Thus, for instance we write ${\sh{\Al}{\element}}^{\sharp} = \Hom_{\qfield}(\sh{\Al}{\element},\qfield)$.
\end{notation}
Let $C=\sh{\Grl}{a}$. The $\completion_\level$-module $\pi^{\level-1}(\ker \ad{a})^{\sharp}=\pi^{\level-1}\Hom_{\completion_{\level}}(\ker \ad{\element},\completion_\level)$ becomes a $\qfield C$-module in a natural way by considering the dual action of the $C$-conjugation by lifts on $\pi^{\level-1}\Al$. Moreover 
\[
\pi^{\level-1}\Hom_{\completion_{\level}}(\ker \ad{\element},\completion_\level)\cong\Hom_{\completion_{\level}}(\pi^{\level-1}\ker \ad{\element},\pi^{\level-1}\completion_\level)
\]
as $\qfield C$-modules via the isomorphism $\pi^{\level-1} \alpha\mapsto \alpha_{\lvert \pi^{\level-1}\ker \ad{\element}}$, and by \cref{rem:change_basis}
\[
\Hom_{\completion_{\level}}(\pi^{\level-1}\ker \ad{\element},\pi^{\level-1}\completion_\level)\cong\Hom_{\qfield}(\redp{\level}(\ker \ad{\element}),\qfield) ={\sh{\Al}{\element}}^{\sharp}.
\]
It then suffices to prove the following:
%
\begin{lem}
\label{lem:coker_ker_dual}
Let $\alginf$ be quadratic with non-degenerate $\mathrm{ad}$-invariant symmetric bilinear form $\form$. Then $\pi^{\level-1}\coker \ad{\element}$ and $\pi^{\level-1}(\ker \ad{a})^{\sharp}$ are isomorphic as $\qfield C$-modules.
\begin{proof}
Consider the dual map of $\ad{\element}$, i.e.\ the map $\ad{\element}^{\sharp}:\Al^{\sharp}\rightarrow \Al$ defined by $f\mapsto f\circ\ad{\element}$. Along the same lines of the proof of \cite[Lemma 8]{jample2012normalforms}, we first prove that $\pi^{\level - 1} \coker \ad{\element}$ and $\pi^{\level - 1} (\ker \ad{a}^{\sharp})^{\sharp}$ are isomorphic as $\qfield C$-modules. The evaluation
\[\alpha_1:\xymatrix@R=3pt{\coker \ad{a}\ar[r]&(\ker \ad{a}^{\sharp})^{\sharp}\\
x+\im \ad{a}\ar@{|->}[r]&(\psi\mapsto \psi(x))}
\]
is an isomorphism of $\completion_\level$-modules and it induces an isomorphism of $\qfield$-vector spaces 
\[
\bar{\alpha}_1: \pi^{\level - 1} \coker \ad{\element}\rightarrow \pi^{\level - 1} (\ker \ad{a}^{\sharp})^{\sharp}.
\]
Moreover $\pi^{\level - 1} (\ker \ad{a}^{\sharp})^{\sharp}$ is a $\qfield C$-module in a natural way by the dual of the $C$-conjugation and one checks that, when $\pi^{\level - 1} (\ker \ad{a}^{\sharp})^{\sharp}$ is equipped with this $\qfield C$-module structure, $\bar{\alpha}_1$ becomes an $\qfield C$-module homomorphism.\par
The second step consists in proving that $\pi^{\level - 1} \ker \ad{a}\cong \pi^{\level - 1} \ker\ad{a}^{\sharp}$ as $\qfield C$-modules. Indeed, if $\form$ is a non-degenerate ad-invariant symmetric bilinear form on $\alginf$, then $\form$ induces a non-degenerate ad-invariant $\completion_\level$-bilinear form $\form_\level$ on $\Al$. This in turn establishes an $\completion_\level$-module isomorphism
\[\alpha_2:\xymatrix@R=3pt{\ker \ad{a}\ar[r]&\ker \ad{a}^{\sharp}\\
x\ar@{|->}[r]&(y\mapsto \form_\level(y,x)),}
\]
and, since $\form_\level$ is ad-invariant, $\alpha_2$ induces an $\qfield C$-module isomorphism
\[
\bar{\alpha}_2: \pi^{\level - 1} \ker \ad{a}\rightarrow \pi^{\level - 1} \ker\ad{a}^{\sharp}.
\]
\end{proof}
\end{lem}
\begin{rem}
\label{rem:explicit_iso}
Under the identification of $\pi^{\level-1}\Al$ with $\algfin$, $\pi^\level \ker\ad{a}$ corresponds to $\sh{\Al}{a}$. Indeed the identification is given by the isomorphism $\varphi: \pi^{\level-1}\Al\rightarrow\algfin$ defined by $\pi^{\level-1} x\mapsto \redp{\level}(x)$. It thus suffices to prove that 
\[\im \varphi_{\lvert \pi^{\level-1} \ker \ad{a}}=\sh{\Al}{\element}.
\]
Let $x\in\centr{\Al}{a}$, and $\bar{x}=\redp{\level}(x)\in\sh{\Al}{a}$. By definition, $\pi^{\level-1} x\in \ker\ad{a}$. Thus $\varphi(\pi^{\level-1} x)=\bar{x}$ and we conclude.
\end{rem}\noindent
Let $\bar{\alpha}_1$ and $\bar{\alpha}_2$ be as in the proof of \cref{lem:coker_ker_dual}. For further usage and convenience of notation we define
\begin{equation}
\label{eq:explicit_iso}
\gamma=\bar{\alpha}_2^{\sharp}\circ \bar{\alpha}_1:\xymatrix@R=3pt{\pi^{\level - 1} \coker\ad{a}\ar[r]&{\sh{\Al}{a}}^{\sharp}\\
\pi^{\level - 1} x + \im \ad{\element}\ar@{|->}[r]&(y\mapsto \form_1(\redp{\level}(x),y)).}
\end{equation}

\subsubsection{Proof of \cref{thm:A}} If $\element$ admits a shadow-preserving lift, then  Lemmata \ref{lem:action_kernel}, \ref{lem:reduction} and \ref{lem:coker_ker_dual}, imply that the $\kersigma$-orbits of elements lying above $\element$ correspond to the elements of $\sh{\Al}{\element}^{\sharp}$. By \cref{cor:linear,lem:coker_ker_dual}, the $\cimage/\kersigma$-action on the set of $\kersigma$-orbits $\sh{\Al}{\element}^{\sharp}$ is the dual of the $\sh{\Grl}{\element}$-conjugation on $\sh{\Al}{\element}$. This proves \cref{thm:A}.


	\subsection{Proof of \cref{thm:B}}
	\label{sec:centr}
Choose $b\in\red^{-1}(a)$ with the same shadow as $\element$ and write $x=b+\pi^{\level} x_c$. 
Replacing $c$ with another element in its same $\sh{\Grl}{a}$-orbit if necessary, we may assume that 
\begin{equation}
\label{eq:thm:centr_1}
\gamma(\varphi_\level(\pi^\level x_c + \im \pi \ad{b}) ) = c.
\end{equation}
Now let $h\in \Redp{\level + 1}^{-1}(\sh{\Grl}{\element})$. As the restriction of the reduction modulo $\primeideal$ to $\centr{\Grlone}{b}$ is surjective onto $\sh{\Grl}{a}$, there is $h'\in\centr{\Grlone}{b}$ such that $h\equiv h'\,\bmod\,\primeideal^\level$, i.e.\ $h = h'\exp(\pi y)$ for some $y\in\Alone$. 
As a result, $h$ acts on $x$ as follows
\begin{align*}
h(b+\pi^{\level} x_c)h^{-1}  & = h'\exp(\pi y)(b+\pi^{\level} x_c)\exp(-\pi y)h'^{-1}\\
					 & = h'(b+\pi^{\level} x_c + \pi [y,b])h'^{-1}\\
					 & = b + h'(\pi^\level x_c + \pi [y,b])h'^{-1}
\end{align*}
It follows that $h$ stabilizes $x$ if and only if $h'$ stabilizes $\pi^\level x_c +\im \pi \ad{b}$ in $\pi^\level\coker\pi \ad{b}$
and, by \cref{eq:thm:centr_1}, this is equivalent to $\Redp{\level + 1}(h')=\Redp{\level + 1}(h)$ stabilizing $c$.
\begin{rem}
\label{thm:centr}
In the notation of \cref{cor:induction}, let $H$ be the kernel of the reduction $\bmod\,\primeideal$ from $\centr{\Grlone}{x}$ to $\sh{\Grlone}{x}$. Then the exponential map establishes a bijection between $\pi \centr{\Alone}{b}$ and $H$.
\end{rem}
\subsection{Proof of  \cref{thm:D}}
\label{sec:nlifts}

Let $e$ be the number of lifts of $\element$ with shadow isomorphic to $\othershadow$ and $f$ be the number of orbits lying above $\element$ whose elements have shadow isomorphic to $\othershadow$. First we show that the cardinality of such orbits only depends on $\element$ and $\othershadow$. Let $b\in\Alone$ be a lift of $\element$ with $\sh{\Grlone}{b}\cong\othershadow$ and let $\hat{\orbit}$ be its $\Grlone$-adjoint orbit. By the orbit-stabilizer theorem and  \cref{thm:centr},
\[
\lvert \hat{\orbit}\rvert = \frac{\lvert \Grlone\rvert}{\lvert \sh{\Grlone}{b}\rvert \, \lvert \centr{\Al}{\element}\rvert}.
\]
So the cardinality of $\hat{\orbit}$ does not depend on the choice of $b$.\par
Let $\orbit$ be the $\Grl$-adjoint orbit of $\element$. All the fibres of the restriction of $\red$ to $\hat{\orbit}$ have the same cardinality, thus $\lvert \hat{\orbit}\rvert /\lvert \orbit\rvert$ is the number of lifts of $\element$ in each $\Grlone$-orbit whose elements have shadow isomorphic to $b$ and that intersects $\red^{-1}(\element)$ non-trivially. It follows that $e = (\lvert \hat{\orbit}\rvert /\lvert \orbit\rvert )f$.\par
Let us expand $\lvert \hat{\orbit}\rvert /\lvert \orbit\rvert$: by the orbit-stabilizer theorem, this is equal to
\[
\frac{\lvert \Grlone\rvert}{\lvert \Grl\rvert}\frac{\lvert\centr{\Grl}{\element}\rvert}{\lvert\centr{\Grlone}{b}\rvert}.
\]
By \cref{thm:centr},
\begin{align*}
\lvert\centr{\Grl}{\element}\rvert	& = \lvert \sh{\Grl}{a}\rvert \cdot \lvert \pi\Al\cap \centr{\Al}{\element} \rvert\\
\lvert\centr{\Grlone}{b}\rvert 	& = \lvert \sh{\Grlone}{b}\rvert \cdot \lvert \centr{\Al}{\element} \rvert,
\end{align*}
and since $\lvert \centr{\Al}{\element}\rvert = \lvert \sh{\Al}{\element} \rvert \cdot \lvert \pi\Al\cap\centr{\Al}{\element}\rvert$, we immediately see that
\[
 \frac{ \lvert \pi\Al\cap \centr{\Al}{\element}\rvert}{\lvert \centr{\Al}{\element} \rvert} = \lvert \sh{\Al}{\element}\rvert^{-1}.
\]
The quantity $ \lvert \sh{\Grl}{a}\rvert / \lvert \sh{\Grlone}{b}\rvert$ is, by \cref{cor:induction}, the size of the $\sh{\Grl}{\element}$-orbit in $\sh{\Alone}{\element}^{\sharp}$ corresponding to $\hat{\orbit}$ by \cref{thm:A}. Therefore, by definition,
\[
\frac{ \lvert \sh{\Grl}{a}\rvert}{\lvert \sh{\Grlone}{b}\rvert} f = \lambda.
\]
By \cref{lem:akov2_2.3,def:dim_shadow}, $\lvert\sh{\Al}{a}\rvert=q^{\dim_{\qfield} \lieshadow}$, while $\frac{\lvert\Grlone\rvert}{\lvert\Grl\rvert}=q^{\dim_{\qfield}\alginf}$ and we conclude.
\subsection{Special linear groups}
When the group scheme in question is a special linear group, \cref{thm:D} may be further refined.  We henceforth set $\alggroup = \SL{\size}{}{}$. According to the notation used so far we define $\alginf = \spl{\size}{(\completion)}{}$, $\algfin = \spl{\size}{(\qfield)}{}$ and, for $\otherlevel\in\N$, $\alg{\otherlevel} = \spl{\size}{(\completion_\otherlevel)}{}$. Let also $\dimalg = n^2 - 1$.\par
The normalized Killing form $\killing$ on $\spl{\size}{(\numfield)}{}$ described in \cite[Section 5]{akov2013representation} is non-degenerate and has integer determinant. If the residue field characteristic of $\completion$  does not divide the determinant of $\killing$, then $\alginf$ is quadratic with non-degenerate ad-invariant bilinear form given by the restriction of $\killing$. This situation happens for all but finitely many places in $\numfield$. From now on $\completion$ is such that $\killing_{\lvert \alginf \times \alginf}$ is non-degenerate. For convenience we shall denote this $\completion$-bilinear form also by $\killing$.\par
\cite[Lemma 2.3]{akov2} tells us that, for special linear groups, the group shadow determines the Lie shadow; we need the following definition in order to precisely state this fact.
\begin{defn}
\label{defn:additive_span}
Let $\level\in\N$. Given a group-shadow $\shadow$, we define
\[\addspan{\shadow}=\linspan(\shadow)\cap\algfin,\]
where $\linspan(\shadow)$ is  the additive span of $\shadow$ when considered as a subset of $\Mat{\size}{(\qfield)}$.
\end{defn}
\begin{lem}[{\cite[Lemma 2.3]{akov2}}]
\label{lem:akov2_2.3}
Assume $q>2$. Let $a\in \Al$ with $\sh{\Grl}{a}=\shadow$, then $\sh{\Al}{a}=\addspan{\shadow}$.
\end{lem}
The next step is to organize shadows by their isomorphism type. We assume for the rest of the section that $q>2$. \Cref{lem:akov2_2.3} legitimates the following definitions:
\begin{defn}
\label{def:allshadows}
For all $\level\in\N$, we choose a set of representatives for the collection of all isomorphism classes of group-shadows of elements in $\alg{\level}$. We denote this set of representatives by  $\allshadows{\alg{\level}}$ and call its members \emph{isomorphism types} of shadows of level $\level$. We also choose a set of representatives for the collection of the isomorphism classes of group shadows of all $\alg{\otherlevel}$ ($\otherlevel\in\N$). We denote this set with
\[
\allshadows{\alginf}
\]
and call its elements isomorphism types of shadows.
Notice that, according to this definition, if $\isoshadow\in\allshadows{\alginf}$ then there are $\otherlevel\in\N$ and $x\in\alg{\otherlevel}$ such that  $\isoshadow=\sh{\grp{\otherlevel}}{x}$, for some $x \in\alg{\level}$.
\end{defn}
\begin{defn}
\label{def:dim_shadow}
Let $\level\in\N$ and  $\isoshadow\in\allshadows{\alg{\level}}$. We define 
\[\dimalg_{\isoshadow}=\dim_{\qfield}\addspan{\isoshadow}.\]
Notice that if $a\in\Al$ and $\sh{\Grl}{a}\cong\isoshadow$, then $\dimalg_{\isoshadow}=\dim_{\qfield} \sh{\Al}{a}$ by \cref{lem:akov2_2.3}.
\end{defn}
\begin{defn}
\label{def:in_shadow}
Let $\level\in\N$, $\isoshadow\in\allshadows{\alg{\level}}$ and $\otherisoshadow\in\allshadows{\alg{\level+1}}$. Let $a\in\Al$ with $\sh{\Grl}{a}\cong\isoshadow$. By \cref{lem:akov2_2.3} we may define
\[
\ninshadow{\isoshadow}{\otherisoshadow}=%
\lvert\lbrace c\in{\sh{\Al}{a}}^{\sharp} \mid \stab_{\isoshadow}{(c)}\cong\otherisoshadow\rbrace\rvert%
= \lvert\lbrace c\in{\addspan{\isoshadow}}^{\sharp} \mid \stab_{\isoshadow}{(c)}\cong\otherisoshadow\rbrace\rvert.
\]
\end{defn}
The last definition does not depend on the choice of $\element$ as the following refined version of \cref{thm:D} explains.
\begin{cor}
\label{prop:n_lifts}
Let $\isoshadow,\otherisoshadow\in\allshadows{\alginf}$. Let  $\level\in\N$ and $\element\in\Al$ with $\sh{\Grl}{\element}\cong\isoshadow$. Assume further that $a\in \Al$ admits a shadow-preserving lift.
Then the number of lifts of $\element$ with shadow isomorphic to $\otherisoshadow$ is equal to
\[
q^{\dimalg-\dimalg_\isoshadow}\,\ninshadow{\isoshadow}{\otherisoshadow}.
\]
\end{cor}
\begin{rem}
\label{rem:nlifts_independent}
The proposition above has the important consequence that the number of lifts of an element of $\element\in\Al$ with shadow isomorphic to $\otherisoshadow$ only depends on (the isomorphism type of) $\sh{\Grl}{\element}$ and on $\otherisoshadow$, not on the choice of $\element$ or on $\level$.
\end{rem}
	\renewcommand{\alggroup}{\mathbf{G}}
	\renewcommand{\alg}[1]{\mathfrak{g}_{#1}}
	\renewcommand{\liering}{\mathfrak{g}}

	\section{Applications to representation zeta functions}
	\label{sec:appl}
	This section contains the proof of \cref{thm:C}. 
We keep the notation of the previous section: so $\alggroup = \SL{\size}{}{}$. It is known  that $\group$ is rigid (i.e.\ its number of continuous complex $i$-dimensional irreducible representations is finite for each $i\in\N$). We say that $\permiss\in\N$ is permissible for $\group$ when $\group^\permiss=\alggroup^\permiss(\completion)$ is potent and saturable (cf.\ \cite[Section 2.1]{akov2013representation} for a definition of potent and saturable groups). By \cite[Proposition 2.3]{akov2013representation} there is $\permiss_0\in\N$ such that $\permiss$ is permissible for $\permiss \geq \permiss_0$. When $\permiss\in\N$ is permissible there is a $\Z_p$-Lie lattice  $\padiclie(\group^\permiss)$ associated with the group $\group^\permiss$. An application of \cref{rem:exp_gl}, \cref{prop:exp_determines_lie_alg} and \cite[Proposition 8.2]{KNVS2011topics} shows that we may identify $\padiclie(\group^\permiss)$ and $\pi^\permiss \alginf$.

\subsection{Kirillov orbit method}
\label{sec:kirillov}
The Kirillov orbit method in  \cite{san2009kirillov} allows to express the representation zeta function as a Poincar\'e series of a matrix of linear forms. We recall the definitions of these objects in a slightly more general setting as this will be useful later.\par
\begin{defn}
\renewcommand{\liering}{\mathfrak{h}}
\renewcommand{\dimalg}{f}
\renewcommand{\halfdim}{e}
Let $\liering$ be a Lie lattice over a ring $\ring$ of $\ring$-rank, say, $\dimalg$. If $\otherbasis = \lbrace b_1,\dots,b_\dimalg\rbrace$ is an $\ring$-basis of $\liering$, for any $b_i,b_j\in \otherbasis$, there are $\lambda_{i,j}^1,\dots,\lambda_{i,j}^\dimalg\in\ring$ such that
\[[b_i,b_j]=\sum_{k=1}^\dimalg \lambda_{i,j}^k b_k.\]
We define the \emph{commutator matrix} of $\liering$ with respect to $\otherbasis$ as
	\begin{equation}
		\label{eq:comm_matrix}
		\Cmatrix{\otherbasis}{\varY}=\left( \sum_{k=1}^\dimalg \lambda_{i,j}^k Y_k\right)_{i,j}\in \Mat{\dimalg}{(\ring[\varY])}
	\end{equation}
with  variables $\varY=(Y_1,\dots,Y_\dimalg)$.\par
\renewcommand{\liering}{\mathfrak{g}}
\renewcommand{\dimalg}{d}
\renewcommand{\halfdim}{h}
\end{defn}
We return now to the previous situation where $\ring = \completion$ and $\mathfrak{h} = \alginf$. Let $\dimalg = \rk_\completion \alginf$. Fix a basis $\basis$ of $\alginf$ and let $\cmatrix = \cmatrix_\basis$. Let now $\level\in\N$ and $\overline{\tuple{w}}\in (\completion/\primeideal^\level)^\dimalg$. Let $\tuple{w}\in \completion^\dimalg$ be a lift of $\overline{\tuple{w}}$. The matrix $\Cmatrix{}{\tuple{w}}$ is an antisymmetric $\dimalg\times\dimalg$ matrix, therefore its elementary divisors may be arranged in $\halfdim=\lfloor \dimalg/2\rfloor$ pairs $(\pi^{a_1},\pi^{a_1}),\dots,(\pi^{a_\halfdim},\pi^{a_\halfdim})$ for $0\leq a_1\leq\dots\leq a_\halfdim\in(\N_0\cup\lbrace \infty\rbrace)$ together with a single extra divisor $\pi^\infty= 0$ if $\dimalg$ is odd. We define
	\begin{align*}
		\nu_{\cmatrix, \level}(\overline{\tuple{w}})			&=(\min \lbrace a_i, \level \rbrace)_{i \in \lbrace 1,\dots, \halfdim \rbrace}.
	\end{align*}
It is easy to see that this definition does not depend on the choice of $\tuple{w}$.
\begin{defn}
\label{eq:N_sets}
Let $\cone{}=(\completion^\dimalg)^*$ and, for $\level\in\N$, $\cone{\level}=((\completion/\primeideal^{\level})^d)^*$. Let $I=\lbrace i_1,\dots,i_\lastindex \rbrace_{<}$ be a (possibly empty) subset of $[\halfdim-1]_0 = \lbrace 0,\dots, \halfdim - 1 \rbrace$ such that $i_1<i_2\dots<i_\lastindex$. We set $i_0=0$ and $i_{\lastindex+1}=\halfdim$ and we write 
	\begin{align*}
		\mu_j		&=i_{j+1}-i_j 							&&\text{for } j\in \lbrace 0,\dots, \lastindex \rbrace; &
		N			&=\sum_{j=1}^\lastindex r_j 				&&\text{for } \mathbf{r}_I=(r_{1},\dots,r_{\lastindex})\in \N^{\lvert I\rvert}.
	\end{align*}
The \emph{Poincar\'e series of $\cmatrix$} is
	\begin{equation*}
		\pseries{\mathcal{R}}{(s)}=\sum_{\stackrel{I\subseteq [\halfdim-1]_0}{I=\lbrace i_1,\dots,i_\lastindex \rbrace_{<}}}\sum_{\mathbf{r}_I\in\N^{\lvert I\rvert}} \lvert \pcoeff{I}{(\cmatrix)}\rvert\,q^{-s\sum_{j=1}^\lastindex r_j(\halfdim-i_j)},
	\end{equation*}
where
	\begin{multline*}
		\pcoeff{I}{(\cmatrix)}= \lbrace \tuple{w}\in \cone{N}\mid \nu_{\cmatrix, N}(\tuple{w})=(\underbrace{0,\dots,0}_{\mu_\lastindex},\underbrace{r_{\lastindex},\dots,r_{\lastindex}}_{\mu_{\lastindex-1}},\\
			\underbrace{r_{\lastindex} + r_{\lastindex-1} ,\dots,r_{\lastindex} + r_{\lastindex-1}, }_{\mu_{\lastindex-2}}%
					\dots%
						,\underbrace{N,\dots,N}_{\mu_0})\in \N_0^\halfdim\rbrace.
	\end{multline*}
If $\basis'$ is another basis for $\alginf$, it is known that $\pseries{\mathcal{R}}{(s)} = \pseries{\mathcal{R}_{\basis'}}{(s)}$, we may therefore define the {\em Poincar\'e series of $\alginf$} as
\[
\mathcal{P}_{\alginf}{(s)} = \pseries{\mathcal{R}}{(s)}.
\]
\end{defn}
The following illustrates the relation between the representation zeta function and the Poincar\'e series.
\begin{prop}[{\cite[Proposition 3.1]{akov2013representation}}]
\label{prop:zeta_poinc}
For all $\permiss$ that are permissible for $\group$
\[
\zetafunc{G^\permiss}=q^{\dimalg\cdot \permiss}\,\mathcal{P}_{\alginf}(s+2).
\]
\end{prop}
\subsection{Poincar\'e series with shadows}
\label{sec:poin_shadows}
We shall rephrase the summation defining the Poincar\'e series so that it fits the language of shadows introduced in \cref{sec:adjoint_orbits}. First of all we relate dimensions of Lie shadows and elementary divisors of the commutator matrix of $\alginf$. If $\otherbasis$ is a basis of an $\ring$-lattice we denote its dual by $\otherbasis^\sharp$.
\begin{prop}
\label{lem:Lie_induction}
Let $\level\in\N$ and let $e\in\Al$. Let $b$ be a shadow-preserving lift of $e$ and $x = b + \pi^\level x_c$ for $x_c\in\Alone$. Let $c = \gamma(\bar{\varphi}_\level(\pi^\level x_c + \pi \im \ad{b} ))$ where  $\gamma$ is as in \cref{eq:explicit_iso} and $\bar{\varphi}_\level$ as in \cref{not:phi}.
	\begin{enumerate}
		\item \label{lem:Lie_induction_a}
		Let $\tuple{e}\in (\completion_\level)^\dimalg$ be the coordinates in $\basis^\sharp$ of $\killing_\level(e,-)$, where $\killing_\level$ is the $\completion_\level$-bilinear form on $\alg{\level}$ induced by reducing $\killing$ modulo $\primeideal^\level$. Then
\[
\dim_{\qfield} \sh{\alg{\level}}{e} = \dimalg - 2 \left\lvert \lbrace a\in\nu_{\cmatrix,\level}(\tuple{e}) \mid a < \level  \rbrace \right\rvert.
\]
		\item \label{lem:Lie_induction_b}
		Let $\mathcal{C}$ be an $\qfield$-basis of $\sh{\Al}{e}$ and let $\tuple{c}$ be the coordinates of $c$ with respect to $\mathcal{C}^\sharp$. Then $\dim_{\qfield} \sh{\Alone}{x} = \dim_{\qfield} \ker \cmatrix_{\mathcal{C}}(\tuple{c})$.
	\end{enumerate}
\begin{proof}
The proof of the first part is a combination of an argument analogous to the one in the proof of \cite[Lemma 3.3]{akov2013representation} with an argument akin to the one on page 148 of \cite{akov2013representation}. The second part is a consequence of the following \cref{lem:stab_lie}.
\end{proof}
\end{prop}
\begin{lem}
\label{lem:stab_lie}
In the notation of \cref{lem:Lie_induction}, let
	\[
	\stab_{\sh{\Al}{e}}(c) = \lbrace y\in \sh{\Al}{e}\mid c([y,z]) = 0 
	\,\forall z\in\sh{\Al}{e}\rbrace.
	\]
Then $\stab_{\sh{\Al}{e}}(c) = \sh{\Al}{x}$
\begin{proof}
Let $\killing_1$ be the $\qfield$-bilinear symmetric form on $\algfin$ induced by $\killing$. By definition $c = \killing_1 (\bar{x}_c, -)$, for $\bar{x}_c\equiv x_c\,\bmod\,\primeideal$. So $y\in\stab_{\sh{\Al}{e}}(c)$ if and only if $\killing_1(\bar{x}_c,[y,z]) = 0$ for all $z\in\sh{\Al}{e}$; i.e.\ if and only if \[
\killing_1([\bar{x}_c,y],z) = 0 \text{ for all $z\in\sh{\Al}{e}$.}
\]
To see that this is equivalent to $y$ being in $\sh{\Alone}{x}$, start by assuming that the latter holds. The shadow determines the Lie shadow and vice-versa (cf.\ \cref{lem:akov2_2.3}), and $b$ is shadow-preserving, so $\sh{\Al}{e} = \sh{\Alone}{b}$. Moreover we may lift $y$ to $\hat{y}\in\centr{\Alone}{x}$, obtaining that $0 = [x,\hat{y}] = [b, \hat{y}] + [\pi^\level x_c, \hat{y}]$.\par 
It follows that $[\pi^\level x_c, \hat{y}]\in\pi\im\ad{b}$ and therefore $\killing_1([\bar{x}_c,y],z) = 0$ for all $z\in\sh{\Al}{e}$, because $\gamma\circ\bar{\varphi}$ is an isomorphism.\par
Conversely assume that $\killing_1([\bar{x}_c,y],z) = 0$ for all $z\in\sh{\Al}{e}$. Choose a lift $\hat{y}$ of $y$ to $\centr{\Alone}{b}$. We have that $\pi^\level[x_c,\hat{y}] = \pi [b,w]$ for some $w\in\Alone$; hence $\hat{y}-\pi w$ centralizes $x$ because $[x, \hat{y}-\pi w ] = [b, \hat{y}] + \pi^\level [x_c, \hat{y}] - \pi [b, w] = 0$.
Since $y \equiv \hat{y} \equiv  \hat{y}-\pi w\,\bmod\primeideal$, $y$ is in the Lie shadow of $x$.
\end{proof}
\end{lem}
We shall need the following notation:
\begin{defn}
\label{def:flags_shadows}
A \emph{decreasing sequence of shadows} is a non-empty set of isomorphism types of shadows 
\[
\lbrace \isoshadow_1,\dots,\isoshadow_\lastindex\rbrace
\]
such that for $0< i < \lastindex$ we have $\dimalg_{\isoshadow_i}>\dimalg_{\isoshadow_{ i + 1}}$ and $\ninshadow{\isoshadow_i}{\isoshadow_{ i + 1}} \neq 0$. The set of all decreasing sequences of shadows is denoted with $\shadowsequences{}$.
\end{defn}

\begin{defn}
\label{def:coefficients}
Let $\shadowsequence= \lbrace \isoshadow_1,\dots,\isoshadow_\lastindex\rbrace\in\shadowsequences{}$ and $\exponents{\shadowsequence}=(r_{\isoshadow_1},\dots, r_{\isoshadow_\lastindex})\in \N^{\shadowsequence}$. Let $N=\sum_{\isoshadow\in \shadowsequence}r_\isoshadow$ and $\cone{N}$ be as in \cref{eq:N_sets}. We define
	\begin{multline*}
		\coeff{\shadowsequence}{\exponents{\shadowsequence}}{\alginf}
		=\left\lbrace x\in\cone{N}\,\left\vert\,\forall \isoshadow_i\in\shadowsequence \,\forall r\in \left(\sum_{j\leq i} r_{\isoshadow_j},\sum_{j\leq i+1} r_{\isoshadow_j}\right] : 		\sh{\grp{\level}}{\reduction{}{\level}(x)}\cong\isoshadow_{i}\right.\right\rbrace.
		\end{multline*}
It is now possible to rewrite the Poincar\'e series: we define
\begin{align*}
\halfdimshadow{\isoshadow}	&=\halfdim-\left\lfloor\frac{1}{2}\dimalg_{\isoshadow}\right\rfloor		&&(\isoshadow \in \allshadows{\alginf}),\\
\shadowsequences{I}		&=\left\lbrace \lbrace \isoshadow_1,\dots,\isoshadow_\lastindex\rbrace\in\shadowsequences{}\,\left\vert\,\halfdimshadow{\isoshadow_j}=i_j\,\forall j\in\lbrace 1,\dots,\lastindex\rbrace\right.\right\rbrace										&&(\text{$I$ as in \cref{eq:N_sets}}).
\end{align*}
Set $\exponents{\shadowsequence}=\mathbf{r}_I$ for all $\shadowsequence\in\shadowsequences{I}$. By \cref{lem:Lie_induction} Part (\ref{lem:Lie_induction_a}), $\pcoeffCard{I}{(\alginf)}=\sum_{\shadowsequence\in\shadowsequences{I}}\lvert \coeff{\shadowsequence}{\exponents{\shadowsequence}}{\alginf}\rvert$. It follows that
\begin{equation}
\label{eq:Poin_shadows}
\mathcal{P}_{\alginf}(s)=\sum_{\shadowsequence\in\shadowsequences{}}\sum_{\exponents{\shadowsequence}\in\N^\shadowsequence} \lvert \coeff{\shadowsequence}{\exponents{\shadowsequence}}{\alginf}\rvert q^{-s\sum_{\isoshadow\in \shadowsequence}r_\isoshadow\cdot\halfdimshadow{\isoshadow}}.
\end{equation}
\end{defn}
\subsubsection{A multiplicative formula for the Poincar\'e series}
We now specialize to $\alggroup = \SL{3}{}{}$. Throughout the rest of this section $\dimalg=8$ and $\halfdim=4$. The normalized Killing form described in \cite[Section 6.1]{akov2013representation}  is non-degenerate for $3\nmid q$. We assume from now on that $3\nmid q$ (beside $q > 2$ as assumed before). We shall now use the results in \cref{sec:nlifts} to give a multiplicative form for the Poincar\'e series of $\spl{3}{(\completion)}$.
\renewcommand{\alginf}{\spl{3}{(\completion)}}
\renewcommand{\algfin}{\spl{3}{(\qfield)}}
\renewcommand{\grpinf}{\SL{3}{(\completion)}{}}
\renewcommand{\grpfin}{\SL{3}{(\qfield)}{}}
\renewcommand{\alg}[1]{\spl{3}{(\completion_{#1})}}
\renewcommand{\grp}[1]{\SL{3}{(\completion_{#1})}{}}
\begin{rem}
\label{rem:central_elements}
Let $\isoshadow\in\allshadows{\alginf}$ and $\lieshadow=\addspan{\isoshadow}$. Let $\basis_\lieshadow$ be an $\completion$-basis for $\lieshadow$ and let $\cmatrix_\lieshadow$ be the commutator matrix of $\lieshadow$ with respect to $\basis_\lieshadow$. Consider the fixed points
\[
\zen{\isoshadow}{\lieshadow^{\sharp}}=\lbrace \omega\in\lieshadow^{\sharp}\mid g.\omega=\omega\,\forall g\in\isoshadow \rbrace\subseteq \lieshadow^{\sharp}
\]
for the action of $\isoshadow$ on $\lieshadow^{\sharp}$. By \cref{prop:n_lifts}, \cref{lem:akov2_2.3} and \cref{lem:Lie_induction} Part (\ref{lem:Lie_induction_b}), $\zen{\isoshadow}{\lieshadow^{\sharp}}$ is the set of elements for which $\cmatrix_\lieshadow$ has rank $0$, and therefore it is an $\qfield$-vector space of dimension $\dimzen{\isoshadow}\in\N_0$, say. This means
\begin{equation*}
\ninshadow{\isoshadow}{\isoshadow}=\lvert\zen{\isoshadow}{\lieshadow^{\sharp}}\rvert=q^{\dimzen{\isoshadow}}.
\end{equation*}
\end{rem}\noindent
\begin{defn}
\label{def:coefficients_f}
Consider $\shadowsequence= \lbrace \isoshadow_1,\dots,\isoshadow_\lastindex\rbrace\in\shadowsequences{}$. Let $\exponents{\shadowsequence}=(r_{\isoshadow_1},\dots, r_{\isoshadow_\lastindex})\in \N^{\shadowsequence}$. Let $\isoshadow_0=\grpfin$ and $\coeff{\shadowsequence}{\exponents{\shadowsequence}}{\alginf}$ be as in \cref{def:coefficients}. We define
\[
\poly{\shadowsequence}{q}=q^{-(\dimalg-\dimalg_{\isoshadow_\lastindex})-\sum_{\isoshadow\in\shadowsequence}\dimzen{\isoshadow}}\cdot\prod_{\isoshadow_{i}\in\shadowsequence} \ninshadow{\isoshadow_{i-1}}{\isoshadow_{i}}.\]
\end{defn}
\begin{lem}
\label{prop:coefficients2}
Let $\shadowsequence$ and $\exponents{\shadowsequence}$ be as in \cref{def:coefficients_f}. Then 
\begin{equation*}
\lvert\coeff{\shadowsequence}{\exponents{\shadowsequence}}{\alginf}\rvert=\poly{\shadowsequence}{q}\cdot\prod_{\isoshadow\in\shadowsequence}\left( q^{\dimalg-\dimalg_{ \isoshadow}+\dimzen{\isoshadow}}\right)^{r_\isoshadow}.
\end{equation*}
\begin{proof}
Since $\spl{3}{(\completion)}{}$ is shadow-preserving we may repeatedly apply \cref{prop:n_lifts}. This together with \cref{rem:central_elements} gives
\begin{equation*}
\lvert\coeff{\shadowsequence}{\exponents{\shadowsequence}}{\alginf}\rvert=\prod_{\isoshadow_{i}\in\shadowsequence} \ninshadow{\isoshadow_{i-1}}{\isoshadow_{i}}\cdot q^{d-\dimalg_{\isoshadow_{i-1}}}\cdot\prod_{\isoshadow\in\shadowsequence}\left( q^{\dimalg-\dimalg_{ \isoshadow}+\dimzen{\isoshadow}}\right)^{r_\isoshadow-1}.
\end{equation*}
The sum $\sum_{\isoshadow_{i}\in\shadowsequence}( \dimalg_{\isoshadow_i}-\dimalg_{\isoshadow_{i-1}})$ is equal to $ \dimalg_{\isoshadow_\lastindex}-\dimalg_{\isoshadow_0}=-(\dimalg-\dimalg_{\isoshadow_\lastindex})$ and we conclude.
\end{proof}
\end{lem}
\Cref{prop:coefficients2} and \cref{eq:Poin_shadows} imply the following:
\begin{equation}
\label{eq:sum_products}
\mathcal{P}_{\alginf}(s)=1 + \sum_{\shadowsequence\in\shadowsequences{}}\poly{\shadowsequence}{q} \cdot\prod_{\isoshadow\in\shadowsequence}\frac{q^{\dimalg-\dimalg_{\isoshadow}+\dimzen{\isoshadow}-s\cdot\halfdimshadow{\isoshadow}}}{1-q^{\dimalg-\dimalg_{\isoshadow}+\dimzen{\isoshadow}-s\cdot\halfdimshadow{\isoshadow}}}.
\end{equation}
Notice that we did not allow the empty set among the shadow sequences, while $I = \emptyset$ was allowed in \cref{eq:N_sets}. This explains the summand $1$ in the equation above.
\begin{rem}
\label{rem:comparison}
We are now able to compare in more detail our methods with the ones in \cite{akov2}. There, the representation zeta function, as said in the introduction, 
is expressed as a sum of some similarity class zeta functions (see \cite[Definition~5.4, Proposition~5.15]{akov2}). These are computed 
recursively in Proposition~6.3, ibid. As a result, the computation boils down to classifying all the group
shadows up to conjugacy and to determining how these behave under lifting. The methods used there to track down the shadow of a lift are mostly ad hoc for each class of shadows.
Our approach, by contrast, essentially uses \cref{thm:D} to perform a simpler, albeit coarser, analogue of such computations. Indeed, 
for two isomorphism types of shadows $\isoshadow$ and  $\otherisoshadow$, the quantities $z_\isoshadow$ and $ \Lambda(\isoshadow, \otherisoshadow)$ 
may be obtained using a commutator matrix of $\addspan{\isoshadow}$. Namely, the former is the number of variables not appearing in that matrix, while 
the latter is given by \cref{lem:stab_lie} and the fact that the Lie shadow and the group shadow determine each other.
\end{rem}
\subsection{Proof of \cref{thm:C}}
Let $\level\in\N$ and $\element\in\alg{\level}$. We say that $\element$ is \emph{regular} if $\dim_{\qfield} \sh{\alg{\level}}{\element}=2$ and that $\element$ is \emph{subregular} if $\dim_{\qfield} \sh{\alg{\level}}{\element}=4$.
\subsubsection{Shadow sequences}
We now determine which are the sequences of shadows that we need to consider to compute \cref{eq:sum_products}.\par
First of all we notice that the set defined in \cref{def:coefficients} does not include $0$ and no other element of $\spl{3}{(\completion_\level)}$ can have shadow equal to that of $0$. It follows that we may exclude decreasing sequences starting with $\SL{3}{(\qfield)}{}$ from those that we need in order to compute \cref{eq:sum_products}.\par
Consider a regular element $a\in \alg{\level}$ on level $\level\in\N$. Its centralizer is abelian\todo{reference for this}, so the action of $\sh{\grp{\level}}{a}$ on ${\sh{\alg{\level}}{a}}^{\sharp}$ is trivial.\par
For what concerns subregular elements we start by considering the situation at level $\level=1$. That is to say, we look at orbits for the action of $\SL{3}{(\qfield)}{}$ on  $\spl{3}{(\qfield)}$. An analysis of the Frobenius rational forms in $\spl{3}{(\qfield)}$ reveals that the possible minimal polynomials of a subregular element are 
\begin{align*}
\minimal{\alpha}	&=(X-\alpha)(X-2\alpha),	&& (\alpha\in\qfield).
\end{align*}
We shall now investigate the isomorphism types of shadows of regular elements. For $i,j\in\lbrace1,2,3\rbrace$, let $e_{i,j}$ be the $3\times 3$ matrix over $\qfield$ with a $1$ in position $(i,j)$ and $0$ everywhere else. It is straightforward to see that if $\alpha\neq 0$ then the corresponding matrix, say $\element\in\spl{3}{(\qfield)}{}$, is semisimple and diagonalizable, so its shadow is isomorphic to $\SubregSem = \sh{\grpfin}{e_{11} + e_{22} - 2 e_{33}}$. All subregular elements that are not semisimple have minimal polynomial $X^2$ i.e.\ they are nilpotent. Let $a\in\algfin$ be such an element, it is an easy computation to show that 
$\sh{\grpfin}{a}\cong \SubregNilp = \sh{\grpfin}{e_{12}}$.\par
In principle we would still need to complete the investigation for shadows appearing only at higher levels; however, since a lift of a subregular element is either regular or shadow-preserving, there cannot be more shadows of subregular elements.\par
We complete $\lbrace \SL{3}{(\qfield)}{}, \SubregSem, \SubregNilp\rbrace$ to a set $\allshadows{\alginf}$ of representatives of isomorphism classes of shadows for all levels. The possible decreasing sequences of shadows that we need to consider are $\lbrace \SubregSem\rbrace$, $\lbrace \SubregNilp\rbrace$, all $\lbrace \otherisoshadow\rbrace$ with $\otherisoshadow\in\allshadows{\alginf}$ such that $\dimalg_\otherisoshadow = 2$ and all   $\lbrace \isoshadow, \otherisoshadow\rbrace$  where $\isoshadow = \SubregSem, \SubregNilp$, $\dimalg_\otherisoshadow = 2$ and $\ninshadow{\isoshadow}{\otherisoshadow} \neq 0$.
\subsubsection{Convention on regular shadows}
\label{sec:reg_convention}
Next we notice that if $\isoshadow$ is the shadow of a regular element ({\em regular shadow} for short), then  $\dimalg_{\isoshadow} = \dimzen{\isoshadow} = 2$ and $\halfdimshadow{\isoshadow} = 3$. This means that for our purposes we need not distinguish regular shadows according to their isomorphism type. More precisely, let $\Reg$ be a symbol distinct from any isomorphism type in $\allshadows{\alginf}$, we define
\begin{gather*}
\begin{align*}
\mathcal{D}'	 	&	= \lbrace \lbrace \SubregSem\rbrace, \lbrace \SubregNilp\rbrace, \lbrace \Reg\rbrace, \lbrace \SubregSem, \Reg\rbrace, \lbrace \SubregNilp, \Reg\rbrace \rbrace;	&
\dimalg'_\isoshadow	&	= \begin{cases}
						2					&\text{ if $\isoshadow = \Reg$}\\
						\dimalg_\isoshadow		 &\text{ if $\isoshadow = \grpfin, \SubregSem, \SubregNilp$;}
					\end{cases}\\
z'_\isoshadow		&	= \begin{cases}
						2					&\text{ if $\isoshadow = \Reg$}\\
						z_\isoshadow		 	&\text{ if $\isoshadow = \grpfin, \SubregSem, \SubregNilp$;}
					\end{cases}	&
\delta'(\isoshadow)	&	= \begin{cases}
						3					&\text{ if $\isoshadow = \Reg$}\\
						\delta(\isoshadow)		 &\text{ if $\isoshadow = \grpfin, \SubregSem, \SubregNilp$;}
					\end{cases}\\
\end{align*}\\
\begin{align*}
\Delta(\isoshadow, \otherisoshadow) &	=\begin{cases}
									\sum_{\stackrel{\otherisoshadow'\in\allshadows{\alginf}}{\dimalg_{\otherisoshadow'}=2}}\ninshadow{\isoshadow}{\otherisoshadow'} &\text{ if $\isoshadow = \grpfin, \SubregSem, \SubregNilp$ and $\otherisoshadow = \Reg$}\\
									\ninshadow{\isoshadow}{\otherisoshadow}		&\text{ if $\isoshadow = \grpfin, \SubregSem, \SubregNilp$ and $\otherisoshadow = \SubregSem, \SubregNilp$;}
								\end{cases}\\
\end{align*}\\
\begin{align*}
\npoly{\lbrace \isoshadow \rbrace}{q}		&	=q^{-(\dimalg-\dimalg'_{\isoshadow})- z'_{\isoshadow}}	\cdot		 \Delta(\grpfin, \isoshadow)	&&\text{ for $\lbrace \isoshadow \rbrace \in	\mathcal{D}'$};\\
\npoly{\lbrace \isoshadow, \otherisoshadow \rbrace}{q}		&	=q^{-(\dimalg-\dimalg'_{\otherisoshadow}) - z'_{\isoshadow} - z'_{\otherisoshadow}} \Delta(\grpfin, \isoshadow) \Delta(\isoshadow, \otherisoshadow)	&&\text{ for $\lbrace \isoshadow, \otherisoshadow \rbrace \in	\mathcal{D}'$.}
\end{align*}
\end{gather*}
Collecting all summands relative to a shadow sequence containing a regular shadow, we may rewrite \cref{eq:sum_products} as
\begin{equation}
\renewcommand{\ninshadow}[2]{\Delta(#1,#2)}
\renewcommand{\dimzen}[1]{z'_{#1}}
\renewcommand{\halfdimshadow}[1]{\delta'(#1)}
\label{eq:sum_products_2}
\mathcal{P}_{\alginf}(s)=1 + \sum_{\shadowsequence\in\mathcal{D}'}\npoly{\shadowsequence}{q} \cdot\prod_{\isoshadow\in\shadowsequence}\frac{q^{\dimalg-\dimalg'_{\isoshadow}+\dimzen{\isoshadow}-s\cdot\halfdimshadow{\isoshadow}}}{1-q^{\dimalg-\dimalg'_{\isoshadow}+\dimzen{\isoshadow}-s\cdot\halfdimshadow{\isoshadow}}}.
\end{equation}
We shall now finish the proof of \cref{thm:C} by computing the ingredients of this last formula.
 \subsubsection{Zeta function}
\renewcommand{\ninshadow}[2]{\Delta(#1,#2)}
\renewcommand{\dimzen}[1]{z'_{#1}}
\renewcommand{\halfdimshadow}[1]{\delta'(#1)}
Let $\element\in\algfin$ be a subregular semisimple element. The orbit of $\element$ has cardinality
\[
\frac{\lvert\SL{3}{(\qfield)}{}\rvert}{\lvert\GL{2}{(\qfield)}{}\rvert}=q^2(q^2+q+1).
\]
Semisimple subregular elements form as many orbits as the possible different minimal polynomials $\minimal{\alpha}$ with $\alpha\neq0$, i.e.\ $q-1$.
Therefore there are 
\begin{equation}
\label{eq:subreg_sem_card}
\ninshadow{\trivialshadow}{\SubregSem}=q^5-q^2
\end{equation}
subregular semisimple elements in total. Moreover, the $\sh{\grpfin}{a}$-action on ${\sh{\algfin}{a}}^{\sharp}$ is the adjoint action of $\GL{2}{(\qfield)}{}$ on $\gl{2}{(\qfield)}$ and as a consequence
\begin{align}
\label{eq:subreg_sem}
\dimalg'_{\SubregSem}&=4, 	&\dimzen{\SubregSem}&=1,& \ninshadow{\SubregSem}{\Reg}&=q\cdot(q^3-1).
\end{align}
\subsubsection{Subregular nilpotent elements}
\label{sec:subreg_nilp}
Choosing the basis  $\otherbasis = \lbrace e_{12}, e_{11} + e_{22} - 2 e_{33}, e_{13}, e_{32}\rbrace$ for $\addspan{\SubregNilp}$, we compute the commutator matrix
\begin{equation*}
\cmatrix_{\otherbasis}(X_0,\dots,X_3)=\begin{pmatrix}
0 & 0 & 0 & 0 \\
0 & 0 & 3 X_{2} & -3 X_{3} \\
0 & -3 X_{2} & 0 & X_{0} \\
0 & 3 X_{3} & - X_{0} & 0
\end{pmatrix}.
\end{equation*}
Let $f = 2$ or $f = 4$. \Cref{lem:Lie_induction} Part \ref{lem:Lie_induction_b} implies that the number elements in $c\in{\sh{\Al}{a}}^{\sharp}$ such that their $\SubregNilp$-stabilizer is isomorphic to $\isoshadow$ with $\dimalg_\isoshadow = f$ is
\[
\lvert\lbrace \tuple{x}\in\qfield^{4} \mid \dim_{\qfield} \ker \cmatrix_{\otherbasis}(\tuple{x}) = f \rbrace\rvert.
\]
So (as we assumend $3\nmid q$) there are $q$ elements of $\addspan{\SubregNilp}^{\sharp}$ on which $\SubregNilp$ acts trivially and $q^4-q$ whose $\SubregNilp$-stabilizer is isomorphic to $\isoshadow$ with $\dimalg_\isoshadow = 4$. This gives us
\begin{align}
\label{eq:subreg_nilp}
\dimalg'_{\SubregNilp}&=4, 	&\dimzen{\SubregNilp}&=1,	&\ninshadow{\SubregNilp}{\Reg}=q\cdot(q^3-1).
\end{align}
The centralizer of a subregular nilpotent element has cardinality $(q-1)q^3$, therefore
\begin{equation}
\label{eq:subreg_nilp_card}
\ninshadow{\trivialshadow}{\SubregNilp}=q^4 + q^3 - q - 1.
\end{equation}
It follows that the number of regular elements at level $1$ is
\begin{equation}
\begin{split}
\label{eq:reg_card}
\ninshadow{\trivialshadow}{\Reg}	&=q^8-1-\ninshadow{\trivialshadow}{\SubregNilp}-\ninshadow{\trivialshadow}{\SubregSem}\\
							&=q \cdot(q - 1) \cdot  (q^6 + q^5 + q^4 - q^2 - 2q - 1).
\end{split}
\end{equation}
\Cref{tab:sl3} gives an overview of the results in equations \cref{eq:subreg_sem,eq:subreg_sem_card,eq:subreg_nilp,eq:subreg_nilp_card,eq:reg_card}.
		\begin{table}
		\caption{Overview for $\SL{3}{(\completion)}{\permiss}$.}
		\label{tab:sl3}
		\begin{tabular}{cccccc}
		\toprule
		$\isoshadow$&$\dimalg'_\isoshadow$&$z'_{\isoshadow}$&$\delta'(\isoshadow)$&$\otherisoshadow$&$\Delta(\isoshadow,\otherisoshadow)$\\
		\midrule
		$\trivialshadow$&8&0&0&$\SubregSem$&$(q^5-q^2)$\\
		&&&&$\SubregNilp$&$(q^4 + q^3 - q - 1)$\\
		&&&&$\Reg$&$q \cdot(q - 1) \cdot  (q^6 + q^5 + q^4 - q^2 - 2q - 1)$\\
		\midrule
		$\SubregSem$&4&1&$2$&$\Reg$&$q\cdot(q^3-1)$\\
		\midrule
		$\SubregNilp$&4&1&$2$&$\Reg$&$q\cdot(q^3-1)$\\
		\midrule
		$\Reg$&2&2&$3$&n.a.&n.a.\\
		\bottomrule
		\end{tabular}
		\end{table}
 With the help of \cref{tab:sl3}, 
applying \cref{eq:sum_products_2} and operating the substitution in \cref{prop:zeta_poinc} we obtain \cref{thm:C}.
\bibliography{../../../Biblio/Database.bbl}

\end{document}